\documentclass[reqno]{amsart}
\pdfoutput=1

\usepackage{amsmath, amssymb}
\usepackage[foot]{amsaddr}
\usepackage{bm}
\usepackage{csquotes}
\usepackage{gensymb}
\usepackage[margin=2.5cm]{geometry}
\usepackage{graphicx}
\usepackage{hyperref}
\usepackage{mathtools}
\usepackage{multirow}
\usepackage{verbatim,amsmath,amsfonts,epsfig,graphics,setspace,amssymb,
float,stmaryrd, multirow}
\usepackage{caption, subcaption}
\usepackage[lined,ruled]{algorithm2e}
\usepackage{enumerate,color,mathtools,epstopdf}
\newcounter{rmrk}[section]

\numberwithin{equation}{section}

\numberwithin{equation}{section}
\newcommand{\set}[1]{\left\{#1\right\}}

\newcommand{\lv}{\left\langle}
\newcommand{\rv}{\right\rangle}

\newcommand{\pspace}{\mathbf{\Lambda}}
\newcommand{\dspace}{\mathbf{\mathcal{D}}}

\newcommand{\model}{\mathcal{M}}
\makeatletter
\renewcommand{\email}[2][]{%
	\ifx\emails\@empty\relax\else{\g@addto@macro\emails{,\space}}\fi%
	\@ifnotempty{#1}{\g@addto@macro\emails{\textrm{(#1)}\space}}%
	\g@addto@macro\emails{#2}%
}
\makeatother
\graphicspath{{media/}}
\title{Goal-oriented adaptive surrogate construction for stochastic inversion}
\author[S. Mattis]{Steven Mattis$^\dagger$}
\author[B. Wohlmuth]{Barbara Wohlmuth}
\address{Chair for Numerical Mathematics, Technical University Munich, Germany}
\address{$^\dagger$ \textnormal{Corresponding author}}
\email{mattis@ma.tum.de, wohlmuth@ma.tum.de}
\date{\today}
\begin{document}
\begin{abstract}
Stochastic inverse problems are generally solved by some form of finite sampling of a space of uncertain parameters.
For computationally expensive models, surrogate response surfaces are often employed to increase the number of samples used in approximating the solution.
The result is generally a trade off in errors where the stochastic error is reduced at the cost of an increase in deterministic/discretization errors in the evaluation of the surrogate.
Such stochastic errors pollute predictions based on the stochastic inverse.
In this work, we formulate a method for adaptively creating a special class of surrogate response surfaces with this stochastic error in mind.
Adjoint techniques are used to enhance the local approximation properties of the surrogate allowing the construction of a higher-level enhanced surrogate.
Using these two levels of surrogates, appropriately derived local error indicators are computed and used to guide  refinement of both levels of the surrogates.
Three types of refinement strategies are presented and combined in an iterative adaptive surrogate construction algorithm.
Numerical examples, including a complex vibroacoustics application, demonstrate how this adaptive strategy allows for accurate predictions under uncertainty for a much smaller computational cost than uniform refinement.
\end{abstract}
\maketitle

%%%%%%%%%%%%%%%%%%%%%%%%%%%%%%%%%%%%%%%%%%%%%%%%%%%%%%%%%%%%%%%%%%%%
%%%%%%%%%%%%%%%%%%%%%%%%%%%%%%%%%%%%%%%%%%%%%%%%%%%%%%%%%%%%%%%%%%%%
\section{Introduction}
\label{sec:Intro}
%%%%%%%%%%%%%%%%%%%%%%%%%%%%%%%%%%%%%%%%%%%%%%%%%%%%%%%%%%%%%%%%%%%%
%%%%%%%%%%%%%%%%%%%%%%%%%%%%%%%%%%%%%%%%%%%%%%%%%%%%%%%%%%%%%%%%%%%%
%The field of Uncertainty Quantification (UQ) has developed into a sophisticated area of research with many contributors over the last few decades.
%Considerable attention and effort has been focused on the development of efficient and accurate UQ methods.
Efficient and accurate methods for uncertainty quantification (UQ) are topics of much interest in the field of computational mathematics and engineering.
While there are many UQ methods that solve a variety of stochastic inverse problems, the most commonly used methods Bayesian inversion \cite{CDS10, Stuart_Bayesian} and ensemble Kalman filtering \cite{Kalman1960, Evensen03}) use Monte Carlo \cite{Robert_Casella_book, Gentle_book} or Markov Chain Monte Carlo  \cite{Gilks_book, Geyer92} sampling techniques to evaluate a Quantity of Interest (QoI) map, which introduces error due to finite sampling.
This error is exacerbated because numerical techniques are used to solve the model instead of solving it exactly.

A common strategy for reducing the effect of finite sampling error is to construct a surrogate to the QoI response surface.
Evaluating the surrogate is done at a greatly reduced computational cost.
%The past several decades has witnessed tremendous advancement in the development and use of surrogates to propagate uncertainties making 
Surrogate modeling is a large topic, so a full review of such techniques is not possible.
%While it is possible to trace back the roots of much of the work involving stochastic finite element approaches \cite{Ghanem_Spanos_book, Ghanem_RedHorse_99} for building surrogates to the seminal papers of Wiener \cite{Wiener_1938} and Cameron and Martin \cite{CM_1947}, 
A widely used class of surrogate approaches involve global polynomial approximations based on stochastic spectral methods \cite{XK_PC_2002, WK_PC_2006, LeMaitre_UQ_04a, LeMaitre_UQ_04b, MNR07, BCW12, BDW13, Prudhomme_Bryant_2015, BPW_2015}.
Another popular approach is using tensor grid and sparse grid stochastic collocation methods for building surrogates \cite{Almeida_Oden_2010, NTW_SparseGrid_2008}, including adaptive methods \cite{ma2009}.
There are also approaches using stochastic optimization to construct global polynomial \cite{li2014} and local \cite{conrad2016} approximations over sequences of distributions adaptively determined from the data.
The surrogate modeling approach considered in this work most closely resembles techniques that exploit derivative information for building piecewise low-order surrogate approximations to improve pointwise accuracy in propagations of uncertainties \cite{MNR07, Estep_Neckels_2006}.

The surrogate response surface is polluted by two sources of error affecting local accuracy \cite{BCW12, BDW13, Prudhomme_Bryant_2015, BPW_2015}. First, there is approximation error introduced by the type of surrogate. Second, there is numerical error in the evaluation of the numerical model used to construct the surrogate.
Both of these are types of discretization errors.
Thus, using a surrogate can represent a trade-off between the reduction in finite sampling error at the expense of an overall increase in the discretization error.
The end result is that our ability to accurately quantify uncertainties by solution of a stochastic inverse problem may be compromised by the use of surrogates unless additional steps are taken to reduce the discretization errors.

Adjoint techniques for finding computable and accurate a posteriori estimates of discretization errors have existed  for decades \cite{Babuska_Rheinboldt_1978, Estep95, Ainsworth_Oden_1997}.
%The derivation of computable and accurate a posteriori estimates of discretization errors based on variational techniques and adjoints dates back several decades \cite{Babuska_Rheinboldt_1978, Estep95, Ainsworth_Oden_1997}. 
Such techniques have served as the basis of the error estimates for polynomial chaos and pseudospectral based surrogates derived in \cite{BCW12}.
Subsequently, in \cite{BDW13, Prudhomme_Bryant_2015}, such error estimates were used as part of a Bayesian inference to quantify uncertainties on parameters to evolutionary partial differential equations where QoI response surfaces were approximated with polynomial chaos techniques and enhanced by the error estimates.

We present an adaptive method for updating the surrogate to accurately make predictions under uncertainty in a stochastic inverse problem setting.
This adaptivity includes increasing the local polynomial order of the approximation, adding more sample points, and increasing the fidelity of the model for certain samples.
The refinement is guided by two levels of surrogate models, with one incorporating adjoint-based a posteriori error estimates to reduce the effect of numerical errors.
The method is designed to simultaneously decrease the effects of both types of discretization errors on the prediction at each iteration.

This paper is organized as follows.
We provide some general notation, terminology, and assumptions used in this work in Section~\ref{sec:notation}, as well as a brief summary of the theory behind stochastic inversion and the various contributions of errors. 
In Section 3, we describe the abstract process of constructing surrogate approximations, identify the various sources of error in the surrogate, and describe the implicit construction of a general piecewise low-order surrogate.
We subsequently provide the conditions relating the exact and surrogate response surfaces for which the approximation solution to the stochastic inverse problem is in fact exact.
%In Section 4, we give a short description of Bayesian inversion and sample-based stochastic inversion
A brief review of adjoint based a posteriori error and derivative estimates along with a list of useful references are provided.
In Section 4, we also describe how we use such error estimates to enhance surrogates by correcting for persistent local biases due to discretization errors.
Such enhanced surrogates are used to derive local error indicators which can be used for local refinement in a variety of ways.
The enhanced surrogates, error indicators, and refinement strategies are combined in an adaptive strategy for surrogate construction.
In Section 5, the method is applied to a number of example problems of varying complexity, including realistic engineering problems.
Conclusions are discussed in Section 6.

\section{Notation, Terminology, and Assumptions}
\label{sec:notation}

%\begin{figure}[htbp]
%\centering
%	\includegraphics[width=0.65\textwidth]{../Figures/Spaces_Surrogate_Construction}
%	\caption{Schematic of mappings between relevant spaces.}\label{Fig:mappings}
%\end{figure}

We present some notation, terminology, and general assumptions for stochastic inversion of a physical system.
Suppose there is a model $\model(u;\lambda)=0$ of the system, where $u$ denotes a vector of state variables determined by the solution of the model for a specified vector of parameters $\lambda$.
These parameters may include coefficients, initial conditions, boundary conditions, source terms, etc.
We assume the space of possible parameters, denoted by $\pspace$, is known, and these parameters explicitly determine the solution to the model.

%A vector of (linear) functionals defines a quantity of interest (QoI) map, denoted by $Q$, and $Q(u(\lambda))$ represents a particular output datum associated to a particular choice of parameters defined by $\lambda$.
A quantity of interest (QoI) map, $Q$, is defined as a vector of linear functions on the model solution, $Q(u(\lambda))$.
%We then write $Q(\lambda):=Q(u(\lambda))$ to emphasize both the dependence of the output data on the model parameters and the fact that in an experimental setting we may be able to control $\lambda$ to observe $Q(\lambda)$ without fully observing $u(\lambda)$. 
Note that $Q$ explicitly depends on $\lambda$.% since $\lambda$ determines the solution of the model.
Let $\dspace :=Q(\pspace)$ denote the space of model QoI.
%This is illustrated in Figure~\ref{Fig:mappings} where a particular solution $u(\lambda)$ in the middle is mapped by the arrow labeled by the QoI map to a model observable on the right denoted by $Q(u(\lambda))$.
%While stochastic inverse problems may be formally posed in infinite-dimensional spaces, when solving them numerically there must be a finite-dimensional discretized problem.
%Therefore, for discretized problems assume $\pspace \subset \mathbb{R}^n$ and $\dspace \subset \mathbb{R}^m$. 
%Let $\mathcal{B}_{\pspace}$ and $\mathcal{B}_{\dspace}$ denote sigma-algebras on $\pspace$ and $\dspace$ respectively.
In general, $\pspace$ and $\dspace$ should be Banach spaces.
Assume that the QoI map defined by $Q$ is piecewise smooth.

Normally, the model is solved using a numerical approximation, resulting in an approximate solution $u_h(\lambda)$ to the model.
Using the approximate solution $u_h(\lambda)$ introduces error into the computation of QoI.
There may be other numerical errors introduced in the QoI being calculated approximately (e.g. using quadrature to approximate an integral quantity or using an iterative solver).
Define $Q_h(\lambda) := Q(u_h(\lambda)))$ as the computed QoI map incorporating these numerical errors and 
$\epsilon_{Q,h}(\lambda) := Q_h(\lambda)-Q(\lambda)$ as the error.
% $$Q_h(\lambda) = Q(\lambda) + \epsilon_{Q,h}(\lambda),$$ where $\epsilon_{Q,h}(\lambda)$ is the numerical error in $Q_h(\lambda)$.
%We may view these errors as a type of perturbation from the exact solution or data. 
%An a priori error analysis often can be used to determine bounds on $\epsilon_{u,h}(\lambda)$ and $\epsilon_{Q,h}(\lambda)$, which define the maximum perturbations to the exact solutions and data we can expect from numerically solving the model.
%This is illustrated by the shaded areas around the exact solutions and QoI data in the middle and right spaces of Figure~\ref{Fig:mappings} which indicates the magnitude of possible perturbations defined by such an a priori error analysis. 
%The solution of a stochastic inverse problem (Bayesian or measure-theoretic) is generally a probability measure $P_{\pspace}$ on $\pspace$.
The numerical solution to the stochastic inverse problem involves a number of (approximate) evaluations of the map $Q$ (i.e. $Q_h$) and gives an approximation to $P_{\pspace}$.

%Let $\dspace$ be the set of possible model output data defined by quantities of interest (QoI) calculated from the solution to the model.
%While stochastic inverse problems may be formally posed in infinite-dimensional spaces, when solving them numerically there must be a finite-dimensional discretized problem.
%Therefore, for discretized problems assume $\pspace \subset \mathbb{R}^n$ and $\dspace \subset \mathbb{R}^m$. 
%Let $\mathcal{B}_{\pspace}$ and $\mathcal{B}_{\dspace}$ denote sigma-algebras on $\pspace$ and $\dspace$ respectively.
%Assume that the QoI map defined by $Q: \pspace \rightarrow \dspace$ is piecewise smooth.

%Knowledge about prior information, observational noise, and probabilities in $\dspace$ are incorporated into a stochastic inversion framework to pose a stochastic inverse problem.

The goal of stochastic inversion could be to approximate the entire probability measure $P_{\pspace}$, but often the goal is to make a prediction under uncertainty.
Suppose $f: \pspace \rightarrow \mathbb{R}$ is an integrable function.
The problem of prediction under uncertainty is to estimate 
\begin{equation} \label{eq:integral}
I = \int_A f(\lambda) dP_{\pspace}
\end{equation}
for some $A \subset \pspace$.
There are a variety of different functions $f$ which could be considered.
For instance, values of $f$ representing characteristic functions could be used to calculate probabilities of sets of interest of parameters, and
values of $f$ representing unobserved quantities of interest could be useful for predictions and decision-making.
In many cases, $I$ can be approximated well without fully resolving $P_{\pspace}$.

In this paper, the stochastic inverse problem is posed as a Bayesian inverse problem; however, the method is generally applicable to other types of stochastic inversion.
Bayesian inversion is an increasingly popular approach.
It has the benefit of generally being well-posed, but is often expensive to implement due to the large number of model evaluations required, so informed adaptivity has the potential to greatly reduce computational costs.
We focus on the applicability of the methods to Bayesian inversion because of its importance in the field and the great potential benefits of such adaptive methods.

A general formation of the Bayesian stochastic inverse problem is described by Stuart \cite{Stuart_Bayesian}.
Suppose that $\pspace$ and $\dspace$ are Banach spaces, and that $Q: \pspace \rightarrow \dspace$ represents the QoI.
In Bayesian inversion, $Q$ is often called the observation operator.
Suppose that $y \in \dspace$ is given data.
The classical inverse problem of finding $\lambda \in \pspace$ such that
\begin{equation}
y = Q(\lambda)
\end{equation}
is typically ill-posed.
However, suppose that the observations $y$ are subject to observational noise.
A more appropriate model is
\begin{equation}
y = Q(\lambda) + \eta,
\end{equation}
where $\eta$ is a mean zero random variable with known statistical properties.

The prior beliefs about $\lambda$ are described in terms of a probability measure $P_0$ on $\pspace$.
Assuming that $\pspace \subset \mathbb{R}^n$ and $\dspace \subset \mathbb{R}^m$,
the probability of $y$ given $\lambda$ has the density
\begin{equation}
\rho(y | \lambda) := \rho(y - Q(\lambda)).
\end{equation}
This is called the data likelihood.
We are interested in the posterior measure $P_{\pspace}$ the probability measure of $\lambda$ given $y$.
Suppose that $\pi_0$ and $\pi_{\pspace}$ are the probability densities associated with $\mu_0$ and $\mu$ respectively.
Using Bayes' formula, 
\begin{equation}
\pi_{\pspace}(\lambda) \propto \rho(y - Q(\lambda)) \pi_0(\lambda).
\end{equation}
There are many methods for sampling from the posterior measure $P_{\pspace}$ using Bayes' formula.
One of the most common is Markov Chain Monte Carlo (MCMC) \cite{Stuart_Bayesian}.
The basic idea of MCMC is to design a Markov chain with the property that a single sequence of output of the chain $\set{\lambda_i}_{i=1}^{\infty}$ is distributed according to $P_{\pspace}$.

%In practice, a finite length chain $\set{\lambda_i}_{i=1}^{M}$ is used.
%There are a variety of numerical methods for forming such chains, notably the Metropolis-Hastings Method.
%A chain representing samples taken from $\mu$ can be used to approximate the probability measure of Voronoi cells in $\pspace$. 
%Given a chain $\set{\lambda_i}_{i=1}^{M}$ and a Voronoi tesselation of $\pspace$, $\set{\mathcal{V}_{i,N}}_{i=1}^N$, the probability measures are given by
%\begin{equation}
%P_{\pspace}(\mathcal{V}_{i,N}) \approx \frac{1}{M} \sum_{j=1}^M \chi_{\mathcal{V}_{i,N}}(\lambda_j).
%\end{equation}
%With this ability to approximate the probability of Voronoi cell events in $\pspace$, we can estimate the local error contributions.
 
%%%%%%%%%%%%%%%%%%%%%%%%%%%%%%%%%%%%%%%%%%%%%%%%%%%%%%%%%%%%%%%%%%%%
\section{Surrogate Models}
\label{sec:Surrogate}
\subsection{Surrogate Models and Error}
A well-known challenge of stochastic inversion is that it often is extremely computationally expensive.
Methods for sampling from and/or describing probability measures such as Markov Chain Monte Carlo (MCMC) and filtering require a \textit{large} number of evaluations of the (approximate) QoI map.
Each evaluation of $Q_h$ requires an evaluation of the model $\mathcal{M}(u; \lambda)$ which is often computationally expensive.
A common approach for reducing the computational cost is by using a surrogate model.
Constructing a surrogate map $Q_s(\lambda)$ often requires using some particular set of samples of $Q_h(\lambda)$ based on a specific type of sampling in $\pspace$, e.g., using a possibly different set of random samples or using deterministic sampling approaches such as sparse grids \cite{Almeida_Oden_2010, NTW_SparseGrid_2008}.
%In surrogate modelling, a smaller number of evaluations of the model $\mathcal{M}(u; \lambda)$ (and possibly auxilliary systems) are performed and used to form an approximation to $Q$ that is extremely computationally cheap to evaluate. 
Let $Q_s(\lambda)$ denote a computationally inexpensive surrogate approximation to $Q(\lambda)$.
The map $Q_s$ is then used for sampling from and/or describing the probability measures of interest.
For the (normal) case where $Q_s$ is constructed from approximate numerical evaluations of $Q$ ($Q_h$) denote the surrogate $Q_{s,h}(\lambda)$.
Let $\epsilon_{s,h}(\lambda)$ denote the error $Q(\lambda) - Q_{s,h}(\lambda)$.
We decompose the error as  $$\epsilon_{s,h}(\lambda) := \epsilon_s(\lambda) + \epsilon_h(\lambda),$$ 
where $\epsilon_s(\lambda)$ is the error in the choice of surrogate due to limited approximation properties of the surrogate, and $\epsilon_h(\lambda)$ is the error in the surrogate from numerical solution of the model.
We now describe a surrogate of piecewise polynomials on Voronoi tessellations.

\subsection{A Piecewise Polynomial Surrogate on an Implicit Voronoi Discretization}
Voronoi tessellations are a convenient way to discretize domains with moderate dimensions, and can be used to define a class of piecewise-defined surrogate models \cite{VPS_Rushdi}.
The space of uncertain parameters $\pspace$ can be discretized by an implicit Voronoi tessellation simply by sampling the space.
Suppose that $\{ \lambda^{(i)}\}_{i=1}^N$ is a finite set of $N$ distinct points in $\pspace$ that we will call ``samples."
Take a metric $d(\cdot, \cdot)$ defined on $\pspace$.
There is a \textbf{Voronoi tessellation} of $\pspace$ denoted by $\{ \mathcal{V}_{i,N} \}_{i=1}^N \subset \pspace$ defined by
$$\mathcal{V}_{i,N} :=  \{ \lambda \in \pspace : d(\lambda^{(i)}, \lambda) \leq  d(\lambda^{(j)}, \lambda),  \forall j=1,2,..., N \}.$$
Each set $\mathcal{V}_{i,N} $ is called \textbf{Voronoi cell}.
Note that two Voronoi cells $\mathcal{V}_{i,N}$ and $\mathcal{V}_{j,N}$ may intersect, but only on a set of measure zero.
This is an implicit tessellation, i.e. the Voronoi cells do not have to be explicitly constructed.
It is generally only necessary to identify which cell $\mathcal{V}_{i,N}$ contains a point $\lambda$ via a nearest neighbor search amongst $\{ \lambda^{(i)}\}_{i=1}^N$.

Suppose that $\set{\lambda_j}_{j=1}^{M}$ is a Markov chain distributed with respect to the probability measure $P_{\pspace}$ (e.g. the output of an MCMC algorithm).
%There are a variety of numerical methods for forming such chains, notably the Metropolis-Hastings Method.
Such a chain can be used to approximate the probability measure of Voronoi cells in $\pspace$.
Given a chain $\set{\lambda_j}_{j=1}^{M}$ and a Voronoi tessellation of $\pspace$, $\set{\mathcal{V}_{i,N}}_{i=1}^N$, the probability measures are given by
\begin{equation}
P_{\pspace}(\mathcal{V}_{i,N}) \approx \frac{1}{M} \sum_{j=1}^M \chi_{\mathcal{V}_{i,N}}(\lambda_j),
\end{equation}
where $\chi_{\mathcal{V}_{i,N}}$ is the characteristic function
 \[
    \chi_{\mathcal{V}_{i,N}}(\lambda)=\left\{
                \begin{array}{lll}
                  1 & \lambda \in \mathcal{V}_{i,N} \\
                  0 & \lambda \not\in \mathcal{V}_{i,N}.
                \end{array}
              \right.
  \]
Such local probability measure estimates can be used for error estimation and adaptivity.

%Let $Q_s(\lambda)$ denote a computationally inexpensive surrogate approximation to $Q(\lambda)$, so $$Q_s(\lambda)=Q(\lambda) + \epsilon_{s}(\lambda),$$ where $\epsilon_s(\lambda)$ is the error in the surrogate.
%We use $Q_s(\lambda)$ to more efficiently map large numbers of samples between $\pspace$ and $\dspace$ as indicated by the dotted arrow mapping between these spaces in Figure~\ref{Fig:mappings}.
%However, constructing $Q_s(\lambda)$ often requires using some particular set of samples of $Q_h(\lambda)$ based on a specific type of sampling in $\pspace$, e.g., using a possibly different set of random samples or using deterministic sampling approaches such as sparse grids \cite{Almeida_Oden_2010, NTW_SparseGrid_2008}.
%We let $Q^{(N)}_{s,h}(\lambda)$ denote this numerically constructed surrogate and $\epsilon_{s,h}(\lambda)$ denote its error. 
%We decompose the error as  $$\epsilon_{s,h}(\lambda) := \epsilon_s(\lambda) + \epsilon_h(\lambda),$$ 
%where $\epsilon_s(\lambda)$ is the error in the choice of surrogate due to limited approximation properties of the surrogate, and $\epsilon_h(\lambda)$ is the error in the surrogate from numerical solution of the model.
%In other words, in practice, the dotted arrow of Figure~\ref{Fig:mappings} is replaced by a mapping that is polluted by multiple sources of error that affect all sets of propagated samples between the spaces. 

Local polynomial approximations of $Q(\lambda)$ on each Voronoi cell could be formed in a variety of ways including interpolating or fitting evaluations of $Q_h$ or using local Taylor expansions.
Interpolants could be useful in relatively low dimensions; however, in higher dimensions a large number of evaluations of $Q_h$ could be required for an accurate approximation, and the polynomials are highly sensitive to the choice of sets of parameters for which the model is solved.
Approximation has some of the same issues as interpolation as well as the additional problem that error is possibly added at points where the model is actually evaluated, i.e. $Q_{s,h}(\lambda_i) \neq Q_h(\lambda_i)$, where $\lambda_i$ are parameters for which the model has been evaluated.
Local Taylor expansions avoid these problems.

To calculate a local Taylor approximation of $Q$ on $\mathcal{V}_{i,N}$, $Q(\lambda)$ and partial derivatives of $Q$ with respect to $\lambda$ must be approximated at some $\lambda_i \in \mathcal{V}_{i,N}$.
An obvious choice of $\lambda_i$ is the generating point of the Voronoi cell $\mathcal{V}_{i,N}$, $\lambda^{(i)}$.
The centroid of $\mathcal{V}_{i,N}$ could be another choice of $\lambda_i$; however, the centroid is not trivial (and possibly quite computationally expensive) to calculate in high dimensions.
Also, if the Voronoi tessellation is refined by adding more generating samples, the centroids would change, and the old centroids may no longer be contained in the same Voronoi cell.
However, for a low dimensional, non-adaptive problem it could be an advantage to use the centroid because of the smaller effective radius.
For an adaptive scheme, approximating the model and derivatives at the generating points of the Voronoi cells is the natural choice.
The simplest local Taylor approximation on Voronoi cells is a piecewise constant approximation.
The piecewise constant surrogate $Q^{(N)}_{s,0}(\lambda)$ is defined as a simple function on a set of Voronoi cells $\set{\mathcal{V}_{i,N}}_{1\leq i\leq N}$ defined by a set of generating samples $\set{\lambda^{(i)}}_{1\leq i\leq N}$, i.e., 
\begin{equation}\label{eq:0surrogate}
	Q^{(N)}_{s,0}(\lambda) = \sum_{1\leq i\leq N}Q_h(\lambda^{(i)})\chi_{\mathcal{V}_{i,N}}(\lambda).
\end{equation}
%where $\chi_{\mathcal{V}_{i,N}}$ is the characteristic function
% \[
%    \chi_{\mathcal{V}_{i,N}}(\lambda)=\left\{
%                \begin{array}{lll}
%                  1 & \lambda \in \mathcal{V}_{i,N} \\
%                  0 & \lambda \not\in \mathcal{V}_{i,N}.
%                \end{array}
%              \right.
%  \]
Constructing $Q^{(N)}_{s,0}$ only requires $Q_h$ to be evaluated at the generating samples $\set{\lambda^{(i)}}_{1\leq i\leq N}$.
The piecewise linear surrogate is defined as
\begin{equation}\label{eq:1surrogate}
	Q^{(N)}_{s,1}(\lambda) = \sum_{1\leq i\leq N} \left[Q_h(\lambda^{(i)}) + 	\nabla_\lambda Q_h(\lambda^{(i)}) (\lambda-\lambda^{(i)})\right] \chi_{\mathcal{V}_{i,N}}(\lambda).
\end{equation}
Constructing $Q^{(N)}_{s,1}$ requires $Q_h$ and $\nabla_\lambda Q_h$ to be evaluated at the points $\set{\lambda^{(i)}}_{1\leq i\leq N}$.
Calculating the Jacobian $\nabla_\lambda Q_h$ is discussed in Section 3.4.
Higher-order piecewise polynomial surrogates may be constructed similarly; however, it is generally computationally prohibitive to explicitly calculate higher-order partial derivatives of $Q_h$ with respect to $\lambda$.
The action of higher order derivative tensors (notably Hessians) on parameters may be computationally viable, but such calculations would have to be performed each time the surrogate is evaluated, so they are not particularly viable for stochastic inversion.

\subsection{Enhancing Surrogates with Error Estimates}
Traditionally, a posteriori error estimates of QoI from differential equation models derived by variational analysis and adjoints were used to guide local $h$- or $p$-adaptivity, i.e., mesh or order refinement, respectively, in the numerical solution to the model (e.g., see \cite{Becker_Rannacher_2001} and the references therein).
%We consider using the adjoint solutions to enhance the surrogate model $Q^{(N)}_{s,h}$. 
%If we can solve the $N_{s1}$ adjoint problems with a higher-order method to produce a set of reliable a posteriori error estimates, $$\set{e_{Q,h}(\lambda^{(k)})}_{1\leq k\leq N_{s1}},$$
Suppose that for each model evaluation $Q_h(\lambda^{(i)})$ used to generate the surrogate, there exists a reliable error estimate $e_{Q,h}(\lambda^{(i)})$.
There is a set of error estimates$\set{e_{Q,h}(\lambda^{(i)})}_{1\leq i\leq N}$ corresponding with the set of samples $\set{\lambda^{(i)}}_{1\leq i\leq N}$.

We can correct for the persistent local bias due to the error $\epsilon_{Q,h}(\lambda^{(k)})$ polluting the evaluation of the surrogate model any $\lambda\in\mathcal{V}_{i,N}$ by enhancing the surrogate model with the error estimates.
We define the enhanced surrogate by
\begin{equation}
\widehat{Q}^{(N)}_{s,h}(\lambda) = Q^{(N)}_{s,h}(\lambda) + \sum_{1\leq i\leq N} e_{Q,h}(\lambda^{(i)}) \chi_{\mathcal{V}_{i,N}}(\lambda).
\label{eq:enhanced_surrogate}
\end{equation}
The error enhanced surrogate has a reduced amount of error due to the numerical solution of the model $\epsilon_h(\lambda)$.
Adjoint problems may be useful for calculating reliable error estimates for QoI and also for calculating derivatives.
\subsection{Adjoint-based a posteriori Error Estimates and Derivatives}
\label{sec:adj}

Suppose the solution to the model $\model(u;\lambda)=0$ is defined by the solution to the finite dimensional parameterized linear system
\begin{equation}\label{eq:pls}
	A(\lambda)\mathbf{u}(\lambda) = \mathbf{b}(\lambda),
\end{equation}
where for each $\lambda\in\pspace\subset\mathbb{R}^m$, $\mathbf{b}(\lambda)\in\mathbb{R}^n$ and $A(\lambda)\in\mathbb{R}^{n\times n}$ is invertible.
Then, for each $\lambda\in\pspace$, there exists a solution $\mathbf{u}(\lambda)\in\mathbb{R}^n$.
Suppose that the QoI map is given by a scalar functional defined by $Q(\lambda)=\lv \mathbf{u}(\lambda), \pmb{\psi} \rv$ where $\pmb{\psi} \in\mathbb{R}^n$ and $\lv \cdot, \cdot \rv$ denotes the inner product.
The adjoint problem to Eq.~\eqref{eq:pls} is 
\begin{equation}\label{eq:pls_adj}
	A(\lambda)^\top\pmb{\phi}(\lambda) = \pmb{\psi},
\end{equation} 
where $\pmb{\phi}(\lambda)$ is the adjoint solution, and $\pmb{\psi}$ is determined by the QoI.
Suppose for a fixed $\lambda\in\pspace$ we numerically solve Eq.~\eqref{eq:pls} to obtain $\mathbf{u}_h(\lambda)\approx \mathbf{u}(\lambda)$ and subsequently compute $Q_h(\lambda)\approx Q(\lambda)$.
%Recall that $$\epsilon_{Q,h}(\lambda):= Q_h(\lambda) - Q(\lambda).$$
%Without the exact value of $Q(\lambda)$, $\epsilon_{Q,h}(\lambda)$ is uncomputable; however, 
The exact error representation is given by
%Using a standard variational analysis and properties of inner products and linear operators, we have
\begin{eqnarray}
	\epsilon_{Q,h}(\lambda) %&=& \lv \mathbf{u}_h(\lambda) - \mathbf{u}(\lambda), \pmb{\psi} \rv \nonumber \\
				 %&=& \lv \mathbf{u}_h(\lambda)- \mathbf{u}(\lambda), A(\lambda)^\top \pmb{\phi}(\lambda) \rv \nonumber \\
				 %&=& \lv A(\lambda)\mathbf{u}_h(\lambda) - A(\lambda)\mathbf{u}(\lambda), \pmb{\phi}(\lambda) \rv \nonumber \\
				 &=& \lv A(\lambda)\mathbf{u}_h(\lambda) - \mathbf{b}(\lambda), \pmb{\phi}(\lambda) \rv. \label{eq:pls_QoI_error}
\end{eqnarray}
%When $\pmb{\phi}(\lambda)$ is given, then Eq.~\eqref{eq:pls_QoI_error} is computable and gives the exact error $\epsilon_{Q,h}(\lambda)$.
Generally, $\pmb{\phi}(\lambda)$ is replaced by an approximation $\pmb{\phi}_h(\lambda)$  in Eq.~\eqref{eq:pls_QoI_error}, giving a computable a posteriori error estimate, which we denote by $e_{Q,h}(\lambda)$.
Typically, we compute $\pmb{\phi}_h(\lambda)$ using a higher order method than used to compute $\mathbf{u}_h(\lambda)$.

%Another use of adjoints is in the efficient computation of derivatives, which we use in constructing the piecewise linear surrogates considered in this work, but in general can play a larger role in sensitivity analysis of the QoI. 
Let $\lambda_i$ denote the $i$th component of the vector $\lambda$ for $1\leq i\leq m$.
Then, differentiating Eq.~\eqref{eq:pls} with respect to $\lambda_i$ and following a similar set of steps, we arrive at
\begin{equation}\label{eq:pls_QoI_deriv}
	\partial_{\lambda_i} Q_h(\lambda) = \lv \partial_{\lambda_i} \mathbf{b}(\lambda) - \left[\partial_{\lambda_i}A(\lambda)\right] \mathbf{u}(\lambda), \pmb{\phi}(\lambda) \rv.
\end{equation}
The partial derivatives of $\mathbf{b}(\lambda)$ and $A(\lambda)$ can often be determined by algorithmic/automatic differentiation, e.g., see \cite{Bartlett_AD_2006}.
Subsequently, this implies that the gradient of the QoI with respect to the parameter $\lambda$, denoted by $\nabla_\lambda Q(\lambda)$, can be approximated by solving both the model and adjoint model exactly once and then computing a finite number of inner products given by Eq.~\eqref{eq:pls_QoI_deriv}.

	This  adjoint-based approach can be applied to most models defined by a linear operator where only a few specific details change.
	For example, when the model is given by a partial differential equation, and a finite element method is used to compute $u_h(\lambda)$, then we generally solve $\phi_h(\lambda)$ either on a refined mesh or using higher order elements to avoid negative effects of Galerkin orthogonality.
	Two comprehensive references on this subject are \cite{Bangerth_book} and \cite{Becker_Rannacher_2001}.
	%The result is an a posteriori error estimate with a similar form to Eq.~\eqref{eq:pls_QoI_error} where the Euclidean inner product is replaced by a duality pairing. % (often an inner product on a Sobolev space).
	%, and the accuracy of the adjoint based error estimate for QoI computed from many different classes of models is well documented \cite{other_papers}.
	%We can also estimate $\nabla_\lambda Q(\lambda)$ in a similar way although Galerkin orthogonality is no longer an issue since no residual appears in Eq.~\eqref{eq:pls_QoI_deriv}, which implies that the same method can be used to compute $\phi_h(\lambda)$ as was used to compute $u_h(\lambda)$ when the model is a differential equation and the goal is to estimate sensitivities of $Q(\lambda)$. 
When the operator defining the model is nonlinear, one must linearize the model operator prior to defining the adjoint problem.

\subsection{General Surrogates}
Our goal is to adaptively form a surrogate $Q^{(N)}_{s,h}$ for use in stochastic inversion.
We want to minimize the number of computationally expensive model and adjoint solves, while still providing an accurate solution to a stochastic inverse problem.
Suppose there are $M$ possible levels to the numerical model. %$\set{l_j}_{j=1}^M$.
Order the levels so that with increasing index $j$ the model fidelity increases.
For example, suppose the forward model $\model(u;\lambda)=0$ represents a 1-D steady-state partial differential equation that we solve numerically using the finite element method.
Suppose we have three different levels of meshes with 10, 100, and 1000 degrees of freedom (DOF) respectively.
Then, $j=1$ is the level associated with the 10 DOF mesh, $j=2$ is the level associated with the 100 DOF mesh, and $j=3$ is the mesh associated with the 1000 DOF mesh.
The error in the QoI due to the numerical solution of the model $\epsilon_h(\lambda)$, should generally decrease with increasing $j$.
Let $Q_{h,j}(\lambda)$ denote the QoI computed using the numerical solution to the model with level $j$ at $\lambda$.
We can construct the piecewise polynomial surrogate using different levels of model evaluations on different Voronoi cells.
Let $l=\set{l_i}_{i=1}^N$ be a set of levels associated with samples $\set{\lambda^{(i)}}_{i=1}^N$, i.e. $1 \leq l_j \leq M$, for $j=1,2,..,N$.
We also can allow different orders of Taylor approximations on different Voronoi cells.
Let $p=\set{p_i}_{i=1}^N$ be a set of local polynomial orders associated with samples $\set{\lambda^{(i)}}_{i=1}^N$.
In practice, we will only allow $p_i=0$ or $p_i=1$ because of the computational cost of using higher local polynomial orders.
We can use these to define a general piecewise polynomial surrogate on Voronoi cells
\begin{equation}\label{eq:general_surrogate}
	Q^{(N)}_{l,p}(\lambda) = \sum_{1\leq i\leq N} \left[Q_{h,l_i}(\lambda^{(i)}) + 	p_i \nabla_\lambda Q_{h,l_i}(\lambda^{(i)}) (\lambda-\lambda^{(i)})\right] \chi_{\mathcal{V}_{i,N}}(\lambda).
\end{equation}
Correspondingly, if there are error estimates, one can have an enhanced general surrogate
\begin{equation}
\widehat{Q}^{(N)}_{l,p}(\lambda) = Q^{(N)}_{l,p}(\lambda) + \sum_{1\leq i\leq N} e_{Q,h}(\lambda^{(i)}) \chi_{\mathcal{V}_{i,N}}(\lambda).
\label{eq:enhanced_general_surrogate}
\end{equation}

%\subsection{Probability Measures on Voronoi Cells}
%Suppose that $\set{\lambda_j}_{j=1}^{M}$ is a Markov chain distributed with respect to the probability measure $P_{\pspace}$ (e.g. the output of an MCMC algorithm).
%%There are a variety of numerical methods for forming such chains, notably the Metropolis-Hastings Method.
%Such a chain can be used to approximate the probability measure of Voronoi cells in $\pspace$. 
%Given a chain $\set{\lambda_j}_{j=1}^{M}$ and a Voronoi tesselation of $\pspace$, $\set{\mathcal{V}_{i,N}}_{i=1}^N$, the probability measures are given by
%\begin{equation}
%P_{\pspace}(\mathcal{V}_{i,N}) \approx \frac{1}{M} \sum_{j=1}^M \chi_{\mathcal{V}_{i,N}}(\lambda_j).
%\end{equation}
%Such local probability measure estimates can be used for error estimation and adaptivity.
%With this ability to approximate the probability of Voronoi cell events in $\pspace$, we can estimate the local error contributions.

%%%%%%%%%%%%%%%%%%%%%%%%%%%%%%%%%%%%%%%%%%%%%%%%%%%%%%%%%%%%%%%%%%%%
%%%%%%%%%%%%%%%%%%%%%%%%%%%%%%%%%%%%%%%%%%%%%%%%%%%%%%%%%%%%%%%%%%%%

\section{Error Estimation and Adaptivity}
\label{sec:Error}
Computed QoI error estimates and derivatives from adjoint methods can be combined  to estimate the error in the solution to the stochastic inverse problem both globally and locally.
Local error estimates can be used to guide a local adaptive scheme to locally improve the accuracy of the surrogate in its relation to the solution of the stochastic inverse problem.
We derive such error estimates, explain how to calculate them numerically, and develop an adaptive refinement scheme based on them.

\subsection{Derivation of Error Indicators}
\label{sec:error_id}
We can use the piecewise polynomial on Voronoi tessellation surrogates on the types of stochastic inverse problems above and in the numerical methods for solving them.
%Stochastic inversion requries a very large number of evaluations of the QoI map $Q(\lambda)$ in order to accurately calculate an approximate probability measure.
%In practice, the exact QoI map cannot generally be evaluated, and a numerical approximation $Q_h(\lambda)$ is evaluated instead.
%This evaluation involves a numerical solve of the model which is often computationally expensive.
It may be computationally prohibitive to evaluate $Q_h$ enough times to get an accurate approximation of $P_{\pspace}$.
The surrogate $Q^{(N)}_{l,p}$ as presented above is extremely cheap to evaluate.
Evaluating $Q^{(N)}_{l,p}(\lambda)$ only requires performing a nearest neighbor search among $\set{\lambda^{(i)}}_{i=1}^N$ and performing some floating point operations.
By replacing the map $Q$ with $Q^{(N)}_{l,p}$ when numerically solving the stochastic inverse problem, it can be solved relatively cheaply.

%There are three main sources of error in the numerical solution of a stochastic inverse problem using the surrogate $Q^{(N)}_{l,p}$.
Let $P_{\pspace}(A)$ be the exact probability of an event $A$ based on the exact map $Q$, and let $P_{\pspace, N,h}(A)$ be the probability of $A$ using the numerical solution to a stochastic inverse problem using the surrogate $Q^{(N)}_{l,p}$ instead of $Q$.
For prediction under uncertainty, the goal is to compute the integral $\int_{A} f d P_{\pspace}$, given a measurable function $f$ and $A$.
Extend $f$ by zero outside of $A$
%\[ \left\{
%\begin{array}{ll}
%f(\lambda) = f(\lambda) & \lambda \in A\\
%f(\lambda) = 0 & \lambda \not \in A.
%\end{array} \right.
%\]
so that $\int_A f dP_{\pspace} = \int_{\pspace} f dP_{\pspace}.$

Let $\mathcal{P}$ be the space of all probability measures on $\pspace$.
We want to find  $Z_\pspace \in \mathcal{P}$ that exactly calculates the integral and also approximates $P_\pspace$ well.
Define the absolute global error with reference to $P_\pspace$  of $Z_\pspace$ as
\begin{equation} \label{eq:global_error}
E_P(Z_\pspace) := \left| \sum_{i=1}^N   \int_{\mathcal{V}_{i,N}} f dP_{\pspace} - \int_{\mathcal{V}_{i,N}} f d Z_{\pspace} \right| + \gamma \left|\sum_{i=1}^N  P_{\pspace}(\mathcal{V}_{i,N}) - Z_{\pspace}(\mathcal{V}_{i,N} ) \right|,
\end{equation}
where $\gamma$ is a Lagrange multiplier.
Any $Z_\pspace$ such that $E_P(Z_\pspace)=0$ is an adequate approximation of $P_\pspace$ for the goal of computing the integral.
Suppose that $\rho_{\pspace}$ is the probability density associated with $P_{\pspace}$ i.e. the Radon-Nikodym derivative of $P_{\pspace}$ with respect to the Lebesgue measure $\mu_{\pspace}$.
%Let $\mathcal{P}_{\rho}$ be the set of Radon-Nikodym derivatives of probability measures in $\mathcal{P}$ with respect to $\mu_{\pspace}$.
%An equivalent problem to finding $Z_\pspace$ such that $E_P(Z_\pspace)=0$, is to find $z_\pspace \in \mathcal{P}_{\rho}$ such that $E_{\rho}(z_\pspace)=0$, 
%where
If $z_\pspace$ is the Radon-Nikodym derivative of $Z_\pspace$ with respect to $\mu_{\pspace}$, then
\begin{equation} \label{eq:global_error2}
E_P(Z_\pspace)=E_{\rho}(z_\pspace):= |E_1(z_\pspace)| + \gamma |E_2(z_\pspace)|,
\end{equation}
where 
\begin{equation}
E_1(z_\pspace) := \sum_{i=1}^N   \int_{\mathcal{V}_{i,N}} \left( f \rho_{\pspace} -  f z_{\pspace} \right) d \mu_{\pspace}
\end{equation}
and
\begin{equation}
E_2(z_\pspace) := \sum_{i=1}^N  \int_{\mathcal{V}_{i,N}} \left(  \rho_{\pspace} -   z_{\pspace} \right) d \mu_{\pspace}.
\end{equation}
%\begin{equation} \label{eq:global_error2}
%E_P(Z_\pspace)=E_{\rho}(z_\pspace):= \left|\underbrace{ \sum_{i=1}^N   \int_{\mathcal{V}_{i,N}} \left( f \rho_{\pspace} -  f z_{\pspace} \right) d \mu_{\pspace} }_{E_1(z_\pspace)} \right|+ \gamma \left|\underbrace{\sum_{i=1}^N  \int_{\mathcal{V}_{i,N}} \left(  \rho_{\pspace} -   z_{\pspace} \right) d \mu_{\pspace} }_{E_2(z_\pspace)}\right| .
%\end{equation}
To balance the error contributions of $E_1$ and $E_2$ in an error reduction algorithm, we want $\left| \frac{\partial E_1}{\partial z_\pspace} \right| \approx \gamma \left| \frac{\partial E_2}{\partial z_\pspace} \right|$.
Let $\pspace^* = \set{\bigcup\limits_{i=1}^N \mathcal{V}_{i,N} | P_\pspace(\mathcal{V}_{i,N}) \neq  Z_\pspace(\mathcal{V}_{i,N} )}$.
If $\pspace^*$ is bounded  then in $\pspace^*$
\begin{equation}
\left| \frac{\partial E_1}{\partial z_\pspace} \right| = \left| \int_{\pspace^*} f d \mu_\pspace \right| \leq \int_{\pspace^*} |f| d \mu_\pspace \text{, and } \left| \frac{\partial E_2}{\partial z_\pspace} \right| = \mu_\pspace(\pspace^*).
\end{equation}
No variation is done in $\pspace \setminus \pspace^*$ since it already has no direct effect on the calculated integral.
Thus, to balance the error contributions
\begin{equation}
\gamma = \frac{1}{\mu_\pspace(\pspace^*)} \int_{\pspace^*} |f| d \mu_\pspace.
\end{equation}
If $(\pspace^* \cap \text{supp}_{\mu_\pspace} (f))$ 
has zero measure, then $\gamma=0$.
$\gamma$ is theoretically undefined if $\pspace^*$ is unbounded, but for most reasonable cases a bounded $\pspace^*$ can be found be refining the Voronoi tessellation.

The global error $E_P$ can be bounded by a sum of local errors $E_i$ on each Voronoi cell using the triangle inequality:
\begin{equation}\label{eq:global_error3}
E_P(Z_\pspace) \leq  \sum_{i=1}^N E_i(P_{\pspace, N,h}) := E(Z_{\pspace}),
\end{equation}
where
\begin{equation}
E_i(P_{\pspace, N,h}) := \left|   \int_{\mathcal{V}_{i,N}} f dP_{\pspace} - \int_{\mathcal{V}_{i,N}} f d Z_{\pspace} \right| + \gamma \left| P_{\pspace}(\mathcal{V}_{i,N}) - Z_{\pspace}(\mathcal{V}_{i,N} ) \right|.
\end{equation}
%\begin{equation} \label{eq:global_error3}
%E_P(Z_\pspace) \leq \sum_{i=1}^N \underbrace{\left( \left|   \int_{\mathcal{V}_{i,N}} f dP_{\pspace} - \int_{\mathcal{V}_{i,N}} f d Z_{\pspace} \right| + \gamma \left| P_{\pspace}(\mathcal{V}_{i,N}) - Z_{\pspace}(\mathcal{V}_{i,N} ) \right| \right)}_{E_i(P_{\pspace, N,h})} = E(Z_\pspace).
%\end{equation}
Suppose the stochastic inverse problem is solved with the surrogate $Q^{(N)}_{l,p}$ and the corresponding enhanced surrogate $\widehat{Q}^{(N)}_{l,p}$.
Denote the probability measure calculated by solving the stochastic inverse problem with $Q_{l,p}$ as $P_{\pspace, N, h}$ and with $\widehat{Q}^{(N)}_{l,p}$ as $\widehat{P}_{\pspace, N, h}$
%Suppose the stochastic inverse problem is solved with the enhanced surrogate $\widehat{Q}^{(N)}_{l,p}$ corresponding with $Q^{(N)}_{l,p}$.
%Denote the probability measure calculated by solving the stochastic inverse problem with $\widehat{Q}_{l,p}$ as $\widehat{P}_{\pspace, N, h}$.
The effect of the deterministic error is smaller in $\widehat{Q}^{(N)}_{l,p}$ which causes less error pollution in $\widehat{P}_{\pspace, N, h}$.
Estimates $\widehat{E}_i$ of $E_i(P_{\pspace, N,h})$ can be computed by replacing the exact probability measure $P_{\pspace}$ with $\widehat{P}_{\pspace, N, h}$ in Equation \ref{eq:global_error3}.
to define the local error indicators
\begin{equation}\label{eq:local_error}
\widehat{E}_i :=\widehat{E}_{int, i} + \widehat{E}_{prob, i},
\end{equation}
where
\begin{equation}
\widehat{E}_{int, i} := \left| \int_{\mathcal{V}_{i,N}} f d\widehat{P}_{\pspace, N, h} - \int_{\mathcal{V}_{i,N}} f d P_{\pspace, N, h} \right|
\end{equation}
and
\begin{equation}
\widehat{E}_{prob, i} = \gamma \left| \widehat{P}_{\pspace,N,h}(\mathcal{V}_{i,N}) - P_{\pspace, N, h}(\mathcal{V}_{i,N} ) \right|.
\end{equation}
%\begin{equation} \label{eq:local_error}
%\widehat{E}_i := \underbrace{ \left| \int_{\mathcal{V}_{i,N}} f d\widehat{P}_{\pspace, N, h} - \int_{\mathcal{V}_{i,N}} f d P_{\pspace, N, h} \right|}_{\widehat{E}_{int, i}} + \underbrace{ \gamma \left| \widehat{P}_{\pspace,N,h}(\mathcal{V}_{i,N}) - P_{\pspace, N, h}(\mathcal{V}_{i,N} ) \right|}_{\widehat{E}_{prob, i}}.
%\end{equation}
In general, these local integrals and the calculation of $\gamma$ must be approximated, and the method of approximation may depend on the model.

\subsection{Approximation of Integrals} 

We approximate integrals in several ways, depending on how computationally expensive $f$ is to evaluate and the measure that we are integrating with respect to.
We look at three cases.
\subsection*{Calculation of $\gamma$ by Emulation}
%In many cases, $f$ may be very cheap to evaluate compared to solving the numerical model, e.g. a characteristic function, a polynomial, algebraic function.
%Suppose this is so and that the solution of the stochastic inverse problem is in the form of a set approximation of $\space$ and computed probability of these sets, e.g. the solution to a measure-theoretic stochastic inverse problem.
%In the case where $f$ is very cheap to evaluate compared to solving the numerical model, e.g. a characteristic function, a polynomial, algebraic function, 
Emulation is the best method to estimate $\gamma$, the weighting factor for $\widehat{E}_{prob, i}$, which is integrated with respect to $\mu_{\pspace}$, the volume measure on $\pspace$.
Let $\set{\lambda_{em}^{(j)}}_{j=1}^{N_{em}}$ be $N_{em}$ uniform (with respect to the volume measure $\mu_{\pspace}$) i.i.d. samples in $\pspace$, and let $\set{\lambda_{em,i}^{(j)}}_{j=1}^{N_{em,i}} = \left( \set{\lambda_{em}^{(j)}}_{j=1}^{N_{em}} \cap \mathcal{V}_{i,N} \right)$ be the $N_{em,i}$ points in $\mathcal{V}_{i,N}$.
We call these ``emulation points" and use them for Monte Carlo integration over $\mathcal{V}_{i,N}$.
%We approximate the integral
%\begin{equation} \label{eq:vol_em}
% \int_{\mathcal{V}_{i,N}} f dP_{\pspace, N,h } = \int_{\mathcal{V}_{i,N} \cap A} f dP_{\pspace, N,h} \approx \frac{P_{\pspace, N, h}(\mathcal{V}_{i,N} \cap A)}{N_{em,i}} \sum_{j=1}^{N_{em,i}} f(\lambda_{em,i}^{(j)}),
%\end{equation}
%and likewise with with $\widehat{P}_{\pspace, N,h }$.
%Thus, we have
%\begin{equation}
%\widehat{E}_{int,i} = \frac{1}{N_{em,i}} \left| \sum_{j=1}^{N_{em,i}}  f(\lambda_{em,i}^{(j)}) \right| \left|\widehat{P}_{\pspace,N,h}(\mathcal{V}_{i,N} \cap A) - P_{\pspace,N,h}(\mathcal{V}_{i,N} \cap A) \right|.
%\end{equation}

In the case where $f$ is very cheap to evaluate compared to solving the numerical model, e.g. a characteristic function, a polynomial, algebraic function, estimate $\gamma$, by Monte Carlo integration:
\begin{equation} \label{eq:gamma}
\gamma \approx \left( \sum_{\substack{i=1 \\ \widehat{P}_{\pspace,N,h}(\mathcal{V}_{i,N}) \not= P_{\pspace,N,h}(\mathcal{V}_{i,N})}}^N \sum_{j=1}^{N_{em,i}} \left| f(\lambda^{(j)}_{em,i}) \right| \right)/ \left( \sum_{\substack{i=1 \\ \widehat{P}_{\pspace,N,h}(\mathcal{V}_{i,N}) \not= P_{\pspace,N,h}(\mathcal{V}_{i,N})}}^N N_{em,i} \right).
\end{equation}
If $f$ is expensive to solve (e.g. requires solving an expensive model) then replace $f(\lambda^{(j)}_{em,i})$ in  (\ref{eq:gamma}) with $f(\lambda_i)$.
We now can calculate $\widehat{E}_{prob,i}$ for each Voronoi cell using its definition in (\ref{eq:local_error}).

\subsection*{Monte Carlo Estimation of $\widehat{E}_{int,i}$ for Cheap Models}
Suppose that $f$ is cheap to evaluate and that the solution of the stochastic inverse problem is a set of points $\set{\lambda_c^{(j)}}_{j=1}^M$ that are distributed according to $P_{\pspace}$, e.g. a chain from an MCMC solution to a Bayesian inverse problem.
Let $\set{\lambda_c^{(j)}}_{j=1}^M$ and $\set{\widehat{\lambda}_c^{(j)}}_{j=1}^{\widehat{M}}$ be the sets distributed according to $P_{\pspace, N,h }$ and $\widehat{P}_{\pspace, N,h }$ respectively.
We approximate the integral 
\begin{equation} \label{eq:mc_int}
\int_{\mathcal{V}_{i,N}} f dP_{\pspace,N,h} \approx \frac{1}{M} \sum_{j=1}^M f(\lambda_c^{(j)})\chi_{\mathcal{V}_{i,N}}(\lambda_c^{(j)}),
\end{equation}
and likewise for $\widehat{P}_{\pspace, N,h }$.
So
\begin{equation}
\widehat{E}_{int,i} = \left|\frac{1}{\widehat{M}} \sum_{j=1}^{\widehat{M}} f(\widehat{\lambda}_c^{(j)})\chi_{\mathcal{V}_{i,N}}(\widehat{\lambda}_c^{(j)}) - \frac{1}{M} \sum_{j=1}^M f(\lambda_c^{(j)})\chi_{\mathcal{V}_{i,N}}(\lambda_c^{(j)}) \right|.
\end{equation}
%$\gamma$ must be calculated using emulated points as shown in Case 1 because of its dependence on the volume measure.
%Then $\widehat{E}_{prob,i}$ can easily be computed.
The estimate of the integral is
\begin{equation}
I_N = \sum_{j=1}^{M} f(\lambda_c^{(j)}).
\end{equation}
$\widehat{I}_N$ can correspondingly be calculated using $\set{\widehat{\lambda}_c^{(j)}}_{j=1}^{\widehat{M}}$.
%\begin{equation}
%\int_A f dP_{\pspace} - \frac{1}{M} \sum_{j=1}^M f(\lambda_j) = 
%\sum_{i=1}^N \int_{\mathcal{V}_i} f dP_{\pspace} -\sum_{i=1}^N \frac{1}{M} \sum_{\substack{j=1 \\ \lambda_j \in \mathcal{V}_i}}^M f(\lambda_j).
%\end{equation}
%where $M_i = \sum_{j=1}^M \chi_{\mathcal{V}_{i,N}}(\lambda_j)$.
%This can be decomposed into a local error identifier
%\begin{equation} \label{eq:mc_er_id}
%E_{I,i} = \sum_{i=1}^N \left(\int_{\mathcal{V}_i} f dP_{\pspace} - \frac{1}{M} \sum_{\substack{j=1 \\ \lambda_j \in \mathcal{V}_i}}^M f(\lambda_j) \right).
%\end{equation}

\subsection*{Estimation of $\widehat{E}_{int,i}$ for Expensive Models}
It is possible that the function $f$ might be computationally expensive to evaluate.
It may involve solving another model or may depend on the solution to the same model that the stochastic inverse problem is based on.
This is common if $f$ represents some kind of model prediction.
In this case, we will approximate $f$ with a simple function approximation on the Voronoi tessellation, i.e.
\begin{equation}
f_N(\lambda) \approx \sum_{i=1}^N f(\lambda^{(i)}) \chi_{\mathcal{V}_{i,N}}(\lambda), \text{   }\forall \lambda \in A.
\end{equation}
The approximation to the integral over a Voronoi cell is
\begin{equation}
\int_{\mathcal{V}_{i,N}} f d P_{\pspace,N,h} = \int_{\mathcal{V}_{i,N} \cap A} f d P_{\pspace,N,h} \approx \int_{\mathcal{V}_{i,N} \cap A} f_N d P_{\pspace,N,h} = f(\lambda^{(i)}) P_{\pspace,N,h}(\mathcal{V}_{i,N} \cap A)=I_{i,N}.
\end{equation}
Likewise, the enhanced integral estimate $\widehat{I}_{i,N}$ can be calculated with $\widehat{P}_{\pspace,N,h}$.
%\begin{equation} \label{eq:int_approx}
%\int_{A} f d P_{\pspace} = \sum_{i=1}^N \int_{\mathcal{V}_{i,N} \cap A} f d P_{\pspace} \approx \sum_{i=1}^N \int_{\mathcal{V}_{i,N} \cap A} f d P_{\pspace, N, h}.
%\end{equation}
%Computationally, we further approximate
%\begin{equation} \label{eq:quadrature}
%\sum_{i=1}^N \int_{\mathcal{V}_{i,N} \cap A} f d P_{\pspace, N, h} \approx\sum_{i=1}^N \int_{\mathcal{V}_{i,N} \cap A} f_N d P_{\pspace, N, h} = \sum_{i=1}^N f(\lambda^{(i)}) P_{\pspace, N, h}(\mathcal{V}_{i,N} \cap A) = I_N.
%\end{equation}
Using the triangle inequality, the local integration error indicator can be expanded:
$$\widehat{E}_{int,i} = \left| \int_{\mathcal{V}_{i,N}} f d \widehat{P}_{\pspace,N,h} - I_{i,N} \right| \leq  \underbrace{\left| \int_{\mathcal{V}_{i,N}} f d \widehat{P}_{\pspace,N,h} - \widehat{I}_{i,N}\right|}_{E_{int,a,i}} + \underbrace{\left|\widehat{I}_{i,N} -  I_{i,N}\right|}_{E_{int,b,i}}$$
%\[
%\begin{array}{lr}
%\int_{\pspace} f(\lambda) d P_{\pspace} - \sum_{i=1}^N f(\lambda^{(i)}) P_{\pspace, N, h}(\mathcal{V}_{i,N}) &  \\
%= \left( \int_{\pspace} f(\lambda) d P_{\pspace} - \sum_{i=1}^{N} f(\lambda^{(i)}) P_{\pspace}(\mathcal{V}_{i,N}) \right) & \text{(I)} \\
%+ \left( \sum_{i=1}^{N} f(\lambda^{(i)}) P_{\pspace}(\mathcal{V}_{i,N}) - \sum_{i=1}^{N} f(\lambda^{(i)}) P_{\pspace,N,h}(\mathcal{V}_{i,N}) \right) & \text{(II)}
%\end{array}
%\]
%\[
%\begin{array}{l}
%\int_{A} f d P_{\pspace} - I_N   \\
%= \underbrace{\left( \int_{A} f d P_{\pspace} - \int_{A} f_N d P_{\pspace} \right)}_{e_{I,a}} 
%+ \underbrace{\left(\int_{A} f_N d P_{\pspace} - I_N \right)}_{e_{I,b}}
%\end{array}
%\]
$E_{int,b,i}$ is easily computable.
By Butler et al. \cite{butler2017}, $E_{int,a,i}$ can be bounded
\begin{equation}
E_{int,a,i} \leq \frac{\widehat{P}_{\pspace,N,h}(\mathcal{V}_{i,N} \cap A)}{\mu_{\pspace}(\mathcal{V}_{i,N} \cap A)} \left| \int_{\mathcal{V}_{i,N} \cap A} f(\lambda) - f(\lambda^{(i)}) d \mu_{\pspace} \right| \leq C_{i,f} \left[2 \sup_{\lambda \in \mathcal{V}_{i,N}} d(\lambda, \lambda^{(i)})\right]
\end{equation}
where  $C_{i,f}$ is a decreasing function of dimension $n$.
If $f$ is Lipschitz continuous with local Lipschitz constant $L_i$, then $C_{i,f}$ is given by
\begin{equation} \label{eq:cif}
C_{i,f} = \frac{L_i \pi^{n/2}\widehat{P}_{\pspace,N,h}(\mathcal{V}_{i,N} \cap A)}{2^n \Gamma(\frac{n}{2}+1)\mu_{\pspace}(\mathcal{V}_{i,N} \cap A)}.
\end{equation}
If the local Lipschitz constant is not known and for (possibly) discontinuous $f$, then $L_i$ in (\ref{eq:cif}) can be estimated by
\begin{equation} \label{eq:Li}
L_i \approx \sup_{\lambda, \gamma \in \mathcal{V}_{i,N} \cap A} \left| \frac{f(\lambda) - f(\gamma)}{\lambda - \gamma} \right|,
\end{equation}
which can approximately be solved by taking $\lambda$ and $\gamma$ from a set of proposal points in $\mathcal{V}_{i,N} \cap A$.
More details can be found in Butler et al. \cite{butler2017}.

\subsection{Types of Refinement}

The calculated local error indicators can be used to refine the surrogate to increase its accuracy and the accuracy for solving the stochastic inverse problem.We consider three types of refinement.

The first type of refinement is \emph{p-refinement}.
%Suppose the surrogate is at refinement stage $k$.
In this type of refinement, the local polynomial order $p_i$ of the surrogate on a Voronoi cell is increased.
This should decrease the local effect of $\epsilon_s(\lambda)$, the error in the approximate QoI map due to the increase in quality of the surrogate model.
Because of the computational cost of calculating higher-order derivatives,
we only consider $p_i=0$ or $p_i=1$, so this refinement can only happen if $p_i=0$.
%Let $\mathcal{I}^p = \set{\mathcal{I}_j}_{i=1}^J$ be a set of $J$ indices corresponding with cells that should be p-refined. 
%Then the gradients $\nabla_{\lambda} Q_{h, l^{(k)}_{\mathcal{I}_j}}(\lambda^{(\mathcal{I}_j)})$ should be calculated and the surrogate should be updated to increase the polynomial order on $\mathcal{V}_{\mathcal{I}_j,N_k}$ for each $\mathcal{I}_j \in \mathcal{I}^p$.
%This new information can be used to to update the surrogate.
p-refinement is outlined in Algorithm \ref{alg:p_refine}.

The second type of refinement is \emph{level-refinement}.
%Suppose the surrogate is at refinement stage $k$.
In this type of refinement, the model level $l_i$ of the surrogate on a Voronoi cell $\mathcal{V}_{i,N_k}$ is increased.
This should decrease the local effect of $\epsilon_h(\lambda)$, the error in the approximate QoI map due to numerical error in solving the model.
A higher-level solve of the model should generally decrease the error in the approximate computations of QoI.
Level refinement is outlined in Algorithm \ref{alg:l_refine}.

The third type of refinement is \emph{h-refinement}.
%Suppose the surrogate is at refinement stage $k$.
In this type of refinement, new samples (generating points for the Voronoi tessellation) are added.
Adding new samples should locally decrease both $\epsilon_s(\lambda)$ and $\epsilon_h(\lambda)$ on the approximate QoI map.
The local Taylor approximations become better because of the decreases radius and the numerical error is extrapolated less.
h-refinement is outlined in Algorithm \ref{alg:h_refine}.

\begin{algorithm}[h]\caption{p-refinement}
\label{alg:p_refine}
Input: Set of indices $\mathcal{I}^p$ of cells to p-refine.

%Solve stochastic inverse problem with $Q^{(N_k)}_{l^{(k)}, p^{(k)}}$ and $\widehat{Q}^{(N_k)}_{l^{(k)}, p^{(k)}}$.
%
%Calculate error estimate $I_E$ and local error indicators $\set{\widehat{E}_i}_{i=1}^{N_k}$.

%Identify set $\mathcal{I}$ of $J$ cells with highest $\widehat{E}_i$ (or all above a tolerance $\epsilon$) and $p_i^{(k)} =0$.
%
\For{$i \in \mathcal{I}_p$}{
Calculate gradient $\nabla_{\lambda}Q_{h,l_i^{(k)}}(\lambda^{(i)})$.

Update $p^{(k+1)}_i = p^{(k)}_i + 1.$

}

\end{algorithm}

\begin{algorithm}[h]\caption{level-refinement}
\label{alg:l_refine}
Input: Set of indices $\mathcal{I}^l$ to level-refine.

%Solve stochastic inverse problem with $Q^{(N_k)}_{l^{(k)}, p^{(k)}}$ and $\widehat{Q}^{(N_k)}_{l^{(k)}, p^{(k)}}$.
%
%Calculate error estimate $I_E$ and local error indicators $\set{\widehat{E}_i}_{i=1}^{N_k}$.

%Identify set $\mathcal{I}$ of $J$ cells with highest $\widehat{E}_i$ (or all above a tolerance $\epsilon$) and $l_i^{(k)} < M$.

\For{$i \in \mathcal{I}^l$}{
Update $l^{(k+1)}_i = l^{(k)}_i + 1.$ \\
Solve model $\mathcal{M}(u; \lambda)$ numerically at $\lambda^{(i)}$ at level $l_i^{(k+1)}$, and calculate corresponding QoIs $Q_{h,l^{(k+1)}_i}(\lambda^{(i)})$ and error estimates $e_{Q,h}(\lambda^{(i)})$. If $p_i^{(k+1)} > 0$, calculate gradient, $\nabla_{\lambda}Q_{h,l_i^{(k+1)}}(\lambda^{(i)})$.

}

\end{algorithm}

\begin{algorithm}[h]\caption{h-refinement}
\label{alg:h_refine}
%Input: Samples $\set{\lambda^{(i)}}_{i=1}^{N_k}$ defining Voronoi tesselation $\set{\mathcal{V}_{i,N_k}}_{i=1}^{N_k}$, levels $l^{(k)}$, polynomial orders $p^{(k)}$, surrogate $Q^{(N_k)}_{l^{(k)}, p^{(k)}}$, enhanced surrogate $\widehat{Q}^{(N_k)}_{l^{(k)}, p^{(k)}}$, and local error indicators $\set{\widehat{E}_i}_{i=1}^{N_k}$.
Input: New samples $\pspace^h$. % $\set{\lambda^j_{new}}_{j=1}^J$ %and corresponding levels $\mathcal{L}_h = \set{l_{new}^j}_{j=1}^J$.

%Solve stochastic inverse problem with $Q^{(N_k)}_{l^{(k)}, p^{(k)}}$ and $\widehat{Q}^{(N_k)}_{l^{(k)}, p^{(k)}}$.
%
%Calculate error estimate $I_E$ and local error indicators $\set{\widehat{E}_i}_{i=1}^{N_k}$.

%Identify $J$ new samples $\set{\lambda_{new}^j}_{j=1}^J$ by (approximately) solving  Eq. (\ref{eq:href_min2}) or Eq. (\ref{eq:href_min}) $J$ times.

Set $N_{k+1} = N_k +J$.

%\For{$i=1,2,.. N_k$}{
%Set $l^{(k+1)}_i = l^{(k)}_i$ and $p^{(k+1)} = p^{(k)}$.
%}

\For{$i=N_k+1,..., N_k + J$}{
Set $\lambda^{(i)} = \lambda_{new}^{i-N_k}$.

Set $l_i^{(k+1)} = l_I^{(k)}$ and $p_i^{(k+1)} = p_I^{(k)}$.

Identify $I$, such that $\lambda^{(i)} \in \mathcal{V}_{I,N_k}$ (via nearest neighbor search).

Solve model $\mathcal{M}(u; \lambda)$ numerically at $\lambda^{(i)}$ at level $l_i^{(k+1)}$, and calculate corresponding QoIs $Q_{h,l^{(k+1)}_i}(\lambda^{(i)})$ and error estimates $e_{Q,h}(\lambda^{(i)})$. If $p_i^{(k+1)} > 0$, calculate gradient $\nabla_{\lambda}Q_{h,l_i^{(k_1)}}(\lambda^{(i)})$.
}

\end{algorithm}

\subsection{Goal-Oriented Adaptive Refinement}
\begin{algorithm}\caption{Goal-Oriented Adaptive Surrogate Construction}\label{alg:ad_samp}
\SetKwInOut{Input}{input}\SetKwInOut{Output}{output}

\Input{Tolerance $\epsilon$ and maximum iterations $its_{max}$.}
\Output{Integral estimate $\widehat{I}_N$.}

Choose initial samples $\set{\lambda^{(i)}}_{i=1}^{N_0}$ defining $\set{\mathcal{V}_{i,N_0}}_{i=1}^{N_0}$, initial polynomial orders $p^{(0)}$, and initial levels $l^{(0)}.$
% = \set{p_i^{(0)}}_{i=1}^{N_0}$, and initial levels $l^{(0)} = \set{l_i^{(0)}}_{i=1}^{N_0}$.

Solve model $\mathcal{M}(u; \lambda)$ numerically at each sample $\lambda^{(i)}$ at level $l_i^{(0)}$ and calculate corresponding QoIs, $Q_{h,l^{(0)}_i}(\lambda^{(i)})$, and error estimates, $e_{Q,h}(\lambda^{(i)})$. If $p_i^{(0)} > 0$, calculate gradient $\nabla_{\lambda}Q_{h,l_i^{(0)}}(\lambda^{(i)})$.

%Set $k=0$ and $I_E > \epsilon$.

Construct  surrogate $Q^{(N_0)}_{l^{(0)}, p^{(0)}}$ and enhanced surrogate $\widehat{Q}^{(N_0)}_{l^{(0)}, p^{(0)}}$, solve the stochastic inverse problem with both. Calculate error estimate $I_E$,  local error indicators $\widehat{E}_i$, and integral estimate $\widehat{I}_N$.

%Solve stochastic inverse problem with $Q^{(0)}_{l^{(0)}, p^{(0)}}$ and $\widehat{Q}^{(0)}_{l^{(0)}, p^{(0)}}$.

%Calculate error estimate $I_E$ and local error indicators $\set{\widehat{E}_i}_{i=1}^{N_0}$.

%Set $k=0$.

\While{$\left| I_E \right| > \epsilon$ \text{ and } $k < its_{max}$}{
%
%Construct  surrogate $Q^{(N_k)}_{l^{(k)}, p^{(k)}}$ and enhanced surrogate $\widehat{Q}^{(N_k)}_{l^{(k)}, p^{(k)}}$ .
%
%Solve stochastic inverse problem with $Q^{(k)}_{l^{(k)}, p^{(k)}}$ and $\widehat{Q}^{(k)}_{l^{(k)}, p^{(k)}}$.
%
%Calculate error estimate $I^{(k)}_E$ and local error indicators $\set{\widehat{E}_i}_{i=1}^{N_k}$.

Add cells with $p_i = 0$ and nonzero probability to $\mathcal{I}^p$ and do   p-refinement (Alg. \ref{alg:p_refine}).

%Set, $\mathcal{I}^l = \emptyset$, $\pspace^h = \emptyset$.

%Do p-refinement for $\mathcal{I}^p$ (Algorithm \ref{alg:p_refine}).

Identify the max. local error, $\widehat{E}_{max}$, and cells to refine, $\mathcal{I}^{ref}$.% = \set{i \in 1 ... N_k | i \not \in \mathcal{I}^p,  \widehat{E}_i > \alpha \widehat{E}_{max}}$.

\For{$i \in \mathcal{I}^{ref}$}{
Calculate $\widehat{E}^l_{\mathcal{J}_i}(i)$,  $\widehat{E}^h_{\mathcal{J}_i}(i)$, and $\lambda^{opt}$.

\eIf{$\widehat{E}^l_{\mathcal{J}_i}(i) \leq \widehat{E}^h_{\mathcal{J}_i}(i)$ and $l_i^{(k)} < l_{max} $}{Add index $i$ to $\mathcal{I}^l$.}{
Add $\lambda^{opt}$ to $\pspace^h$.}
}

Level-refine (Alg. \ref{alg:l_refine}) for $\mathcal{I}^l$, h-refine (Alg. \ref{alg:h_refine}) for $\pspace^h$, and update iteration number $k = k +1$

Construct  surrogate $Q^{(N_k)}_{l^{(k)}, p^{(k)}}$ and enhanced surrogate $\widehat{Q}^{(N_k)}_{l^{(k)}, p^{(k)}}$ and solve the stochastic inverse problem with both. Calculate error estimate $I_E$,  local error indicators $\widehat{E}_i$, and integral estimate $\widehat{I}_N$.
%
%Solve stochastic inverse problem with $Q^{(k)}_{l^{(k)}, p^{(k)}}$ and $\widehat{Q}^{(k)}_{l^{(k)}, p^{(k)}}$.

%Calculate error estimate $I^{(k)}_E$ and local error indicators $\set{\widehat{E}_i}_{i=1}^{N_k}$.
}
\end{algorithm}

The computed probabilities, enhanced probabilities, and local error indicators can be used to guide adaptive refinement for the goal of accurately calculating the integral.
Consider an initial discretization via an implicit Voronoi tessellation $\set{\mathcal{V}_{i,N_0}}_{i=1}^{N_0}$ of $\pspace$ by $N_0$ points $\set{\lambda^{(i)}}_{i=1}^{N_0}$.
%The subscript refers to the refinement stage
Solve the model with these input parameters, at the lowest level and calculate error estimates (and possibly derivatives).
The initial levels are $l^{(0)} = \set{1}_{i=1}^{N_0}$.
Choose the polynomial order $p_0$ (0 or 1) with which to define the surrogate depending on whether or not derivatives are known.
The initial polynomial orders are $p^{(0)} = \set{p_0}_{i=1}^{N_0}$.
Thus, we can construct the initial surrogate $Q^{(N_0)}_{l^{(0)}, p^{(0)}}(\lambda)$ and enhanced surrogate $\widehat{Q}^{(N_0)}_{l^{(0)}, p^{(0)}}(\lambda)$ on $\pspace$.
Starting with the initial surrogate and enhanced surrogate, an iterative procedure is performed for goal-oriented adaptive refinement.
We denote each iteration with an index $k$.

If an adjoint approach is being used to calculate error estimates and to enhance the surrogate, then derivative information is computationally cheap to obtain as discussed in Section \ref{sec:adj}.
A locally piecewise linear surrogate has much more accuracy compared to the piecewise constant surrogate, so if derivative information is known, p-refinement (from p=0 to p=1) should be performed for all cells with nonzero probability at each iteration.
%Set 
%\begin{equation} \label{eq:Ip}
%\mathcal{I}^p = \set{i \in 1 ... N_k | p_i^{(k)} = 0 \text{ and } \widehat{P}_{\pspace, N,h}(\mathcal{V}_{i,N_k})  \neq 0 \text{ or } P_{\pspace, N,h}(\mathcal{V}_{i,N_k})  \neq 0 },
%\end{equation}
%and perfrom p-refinement on those cells before doing any other type of refinement on them. 

The local error indicators $\set{E_i}_{i=1}^{N_k}$ are used to guide h- and level-refinement.
In practice, the exact error indicators $\set{E_i}_{i=1}^{N_k}$ are not known, so we use  the approximations $\set{\widehat{E}_i}_{i=1}^{N_k}$ as described in Section \ref{sec:error_id}.
%Doing this, the numerical implentations of p-refinement and level-refinement are straightforward.
%The $J$ samples with the greatest error indicator or have error indicators above some tolerance that have the ability to be refined (are not at the highest level or highest polynomial order) are refined.
We want to perform h- or level-refinement on the Voronoi cells with the highest local error.
The maximum local error indicator is $\widehat{E}_{max} = \max\limits_{i=1,..,N_k} \widehat{E}_i$.
Given a parameter $\alpha \leq 1$, we want to refine all cells $i$ such that $\widehat{E}_i > \alpha \widehat{E}_{max}$.
Denote the set of these indices as $\mathcal{I}^{ref}$.
There are alternative techniques for marking cells for refinement (such as the mean strategy or D\"{o}rfler strategy which are common in adaptive finite element methods \cite{demkowicz2006}).
The type of refinement for each cell $i \in \mathcal{I}^{ref}$ must be determined.

Determining between h- and level-refinement for a sample should be done by determining which type (potentially) reduces the error in a neighborhood of the cell the most.
Let $\mathcal{J}_i \subset \set{1,2, .., N_k} $ be a set of indices of cells in a neighborhood of cell $i$.
These are cells whose local error indicators would likely be changed if refinement is done for cell $i$.
Obvious choices of $\mathcal{J}_i$ are $i$ and its direct neighbors or cells within some distance of $\lambda^{(i)}$.
Let $\widehat{E}_j^l(i)$ be defined as the local error indicator for cell $j$ after level-refinement is done for cell $i$.
The total local error sum of over $\mathcal{J}_i$ under level-refinement is 
\begin{equation}
 \widehat{E}^l_{\mathcal{J}_i}(i) =  \sum_{j \in \mathcal{J}_i} \left(\widehat{E}_j^l(j)\right).
\end{equation}

Let $\widehat{E}_j^h(\lambda)$ be defined as the local error indicator for cell $j$ after a sample $\lambda$ is added to the discretization via h-refinement.
%Let $\pspace_{local,i} \subset \pspace$, with $\mathcal{V}_{i, N_k} \subset \pspace_{local,i} $ be a region in $\pspace$ such that if if a point contained in it is added to the surrogate $\mathcal{V}_{i, N_k}$ will change, i.e. the set of points that will h-refine $\mathcal{V}_{i, N_k}$. 
The total local error sum over $\mathcal{J}_i$ under h-refinement with $\lambda$ is 
\begin{equation}
 \widehat{E}^h_{\mathcal{J}_i}(\lambda) =  \sum\limits_{j \in \mathcal{J}_i} \left(\widehat{E}_j^h(\lambda)\right) + \widehat{E}_{N_k+1}^h(\lambda),
\end{equation}
and the optimal sample to add is
\begin{equation} \label{eq:opt_h_point}
\lambda^{opt} = \arg \min_{\lambda \in \pspace} \widehat{E}^h_{\mathcal{J}_i}(\lambda).
\end{equation}
If $\widehat{E}^l_{\mathcal{J}_i}(i) \leq \widehat{E}^h_{\mathcal{J}_i}(\lambda^{opt})$ level-refinement should be done for cell $i$, and otherwise h-refinement should be done by adding $\lambda^{opt}$ to the surrogate.

In practice, these optimization problems are unfeasible to solve directly, but simple approximations can be used to estimate $\widehat{E}^l_{\mathcal{J}_i}(i)$, $\lambda^{opt}$, and $\widehat{E}^h_{\mathcal{J}_i}(\lambda^{opt})$.
The effect of level-refinement can be approximated by locally replacing $P_{\pspace, N,h}$ by $\widehat{P}_{\pspace, N,h}$, rescaling (assuring it integrates to the original value) $P_{\pspace, N,h}$ on the rest of the neighborhood, and recalculating the local error indicators.
The effect of h-refinement can approximated by adding a reasonable number of proposal points to the Voronoi tessellation and recalculating the local error indicators.
The proposal with the smallest local error indicator is the best proposal point for h-refinement.
%The process is outlined in Algorithm \ref{alg:hp-approx}.

After all of the cells that will undergo refinement and the type of refinements are identified, the model is solved correspondingly. For p-refinement, derivatives are calculated.
For level-refinement, the model and adjoints are solved with the higher level model and corresponding QoIs, error estimates, and derivatives are calculated for the parameters corresponding with the generating point of the Voronoi cell.
For h-refinement the model and adjoints are solved and corresponding QoIs, error estimates, and derivatives are calculated for the new parameter.
This new information is used to update the surrogate and enhanced surrogate to stage $k+1$.
The stochastic inverse problem is solved, local error indicators are calculated, and new iteration of refinements are performed.
The process terminates when a stopping criterion is met.
A reasonable stopping criterion is when the approximation of the integral has not varied within some tolerance for several steps.
If computational resources are limited the process should also be stopped after some computational budget is met.
The adaptive method, including surrogate construction, was implemented using the open-source Python package BET \cite{pyBET}.

\section{Numerical Results}
\label{sec:Numerics}
We have applied our adaptive scheme for goal-oriented surrogate construction to a variety of example problems.
The first example is a simple 1D second order PDE system with two uncertain parameters.
This simple low-dimensional problem helps illustrate the algorithm and results can be displayed visually.
This problem is cheap to solve and has a analytical solution; however, the relationships between parameters and QoIs are highly nonlinear so it useful for illustrating the benefits of the method.
The second example is an 2-D elliptic boundary value problem with a complicated conductivity field parameterized by eight coefficients.
The third is a nonlinear system of ordinary differential equations where two initial conditions and four coefficients are uncertain.
The fourth and final example is a complicated engineering problem from vibroacoustics involving the deformation of a violin bridge.
\subsection{1D Elliptic PDE}
The first example is a one-dimensional elliptic boundary value problem with homogeneous Dirichlet boundary conditions
\begin{eqnarray}
-\lambda_1 v''(x) &=& \exp(\lambda_2 x), \text{ for } x\in (0,1) \\
v(0)&=&v(1)=0, \nonumber
\end{eqnarray}
where $\lambda =[\lambda_1, \lambda_2]$ are the uncertain parameters.
Let $u(\lambda, x) = v(x)$ be the solution to the problem with given parameters $\lambda$.
Suppose that that QoI map is $Q(\lambda) = [Q_1(\lambda), Q_2(\lambda)]$ where $Q_1(\lambda) = \int_{0.1}^{0.4} u(\lambda, x) dx$ and $Q_2(\lambda) = \int_{0.6}^{0.9} u(\lambda, x) dx$.
Let $\psi_1$ and $\psi_2$ be the Riesz representors of $Q_1$ and $Q_2$ in $L^2(0,1)$.
For reference, $u(\lambda,x)$, $Q_1(\lambda)$ and $Q_2(\lambda)$ can be expressed exactly by 
\begin{equation}
u(\lambda,x) = {\frac {-{{\rm e}^{\lambda_2\,x}}-x+x{{\rm e}^{\lambda_2}}+1}{\lambda_1{\lambda_2}^{2}},
}
\end{equation}
\begin{equation}
Q_1(\lambda) = {\frac { 0.075\,{{\rm e}^{\lambda_2}}\lambda_2-{{\rm e}^{ 0.4\,\lambda_2}}+{
{\rm e}^{ 0.1\,\lambda_2}}+ 0.225\,\lambda_2}{\lambda_1{\lambda_2}^{3}}},
\end{equation}
and
\begin{equation}
Q_2(\lambda) = {\frac { 0.225\,{{\rm e}^{\lambda_2}}\lambda_2-{{\rm e}^{ 0.9\,\lambda_2}}+{
{\rm e}^{ 0.6\,\lambda_2}}+ 0.075\,\lambda_2}{\lambda_1{\lambda_2}^{3}}}.
\end{equation}
These exact solutions will be used as a reference for comparison with results using a numerical method.

The problem can easily be discretized with a standard centered finite difference approximation with uniform spacing $h$ forming a linear system
\begin{equation}
A_h(\lambda) \mathbf{u}_h(\lambda) = \mathbf{b}_h(\lambda)
\end{equation}
for given parameters $\lambda$.
This system can be efficiently solved directly using a sparse solver.
$\psi_{1}$ and $\psi_{2}$ can be discretized as $\psi_{h,1}$ and $\psi_{h,2}$ and can be used to evaluate the approximate QoI map $Q_h(\lambda) = [Q_{h,1}(\lambda), Q_{h,2}(\lambda)]$, where $Q_{h,1}(\lambda) = \lv \mathbf{u}_h(\lambda), \psi_{h,1} \rv$ and $Q_{h,2}(\lambda) = \lv \mathbf{u}_h(\lambda), \psi_{h,2} \rv$.
We solve the adjoint problems on a mesh that is twice as fine
\begin{equation}
A_{h/2}(\lambda)^T \boldsymbol{\phi}_{h/2,1}(\lambda) = \boldsymbol{\psi}_{h/2,1}
\end{equation}
\begin{equation}
A_{h/2}(\lambda)^T \boldsymbol{\phi}_{h/2,2}(\lambda) = \boldsymbol{\psi}_{h/2,2}.
\end{equation}
The adjoint solutions $\boldsymbol{\phi}_{h/2,1}(\lambda)$ and $\boldsymbol{\phi}_{h/2,2}(\lambda)$ are used to calculate $\epsilon_{Q,h}(\lambda)$ and $\partial_{\lambda} Q_h(\lambda)$ as shown in (\ref{eq:pls_QoI_error}) and (\ref{eq:pls_QoI_error}) respectively.
There are five levels of model resolution corresponding with meshes with $h=0.2$, $h=0.1$, $h=0.05$, $h=0.025$, and $h=0.0125$ respectively.
This provides all of the ingredients necessary to form regular and enhanced piecewise constant, piecewise linear, and general surrogates $Q^{(N)}_{l,p}$.

\begin{figure}
\centering
\includegraphics[width=0.8\textwidth]{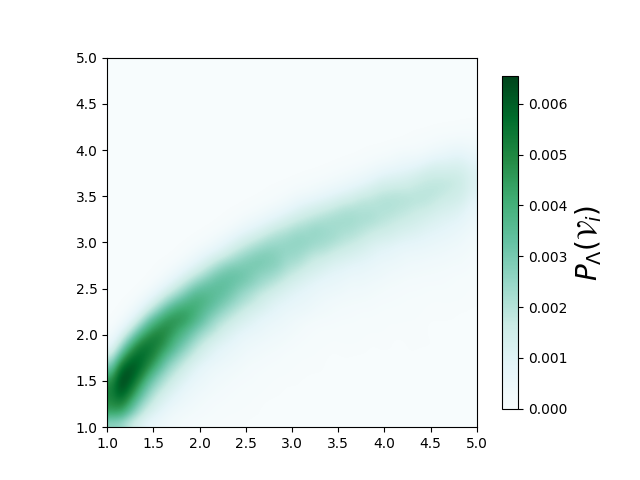}
\caption{Reference posterior distribution.}
\label{fig:true_post}
\end{figure}

The stochastic inverse problem is a Bayesian inverse problem with the forward map $Q: \pspace \rightarrow \mathbb{R}^2$, where $\pspace = [1,5]^2$ and $Q$ is defined above.
The data is $y = [0.22, 0.15]$.
We assume a uniform prior on $\pspace$ and mean-zero Gaussian noise $\eta \backsim \mathcal{N}([0,0], [0.0025, 0.0025])$.
The function $f$ that we are interested in integrating with respect to the posterior is $f(\lambda) = \frac{dv}{dx}\vert_{0.83}$.
A standard Metropolis-Hastings MCMC method is used for solving the stochastic inverse problem.
A reference calculation of the posterior $P_{\pspace}$ is calculated using the exact evaluation of $Q$ and an MCMC solution with $10^7$ samples and is shown in Figure \ref{fig:true_post}.
Using the true posterior, the true value of the integral is $\int_{\pspace} f P_{\pspace} = -0.60178$ calculated using Monte Carlo integration.
%The true posterior is shown in Figure \ref{fig:true_post}.

For comparison with the adaptive algorithm, the integral was calculated using the posterior calculated using piecewise constant and piecewise linear surrogates with uniformly distributed generating samples.
10, 100, 1000, and 10000 uniform i.i.d. generating points were used and the model was solved with all five levels of the mesh.
Tables \ref{tab:ex1_const} and \ref{tab:ex1_lin} show the average absolute errors in the calculation of the integral using piecewise constants and piecewise linears respectively with the different numbers of generating points and model levels.
We see that as the mesh level is increased the calculations usually become more accurate because $\epsilon_h$ is decreasing.
However, when there are a small number of generating samples, the effect of the surrogate error $\epsilon_s$ dominates and the convergence slows down or end completely.
There is also the tendency for the error to decrease as the number of generating samples increases, because $\epsilon_s$ is decreasing.
However, at the low levels there is more deterministic model error and the convergence slows or stops.
This is because $\epsilon_h$ is polluting the calculation, and the effect cannot be removed without reducing the error in the solution of the QoIs.
The similarity in errors between the piecewise constant and piecewise linear surrogates is explained in an analogous way to why the trapezoidal rule and midpoint rules for quadrature in one-dimension are the same order.
We see this tendency go away in higher dimensional problems.
This illustrates that for the solution of the integral to be truly accurate, the surrogate needs to reduce the effects of both of these types of error.
The adaptive strategy is designed to reduce both of these effects efficiently.

\begin{figure}
\centering
\includegraphics[width=0.8\textwidth]{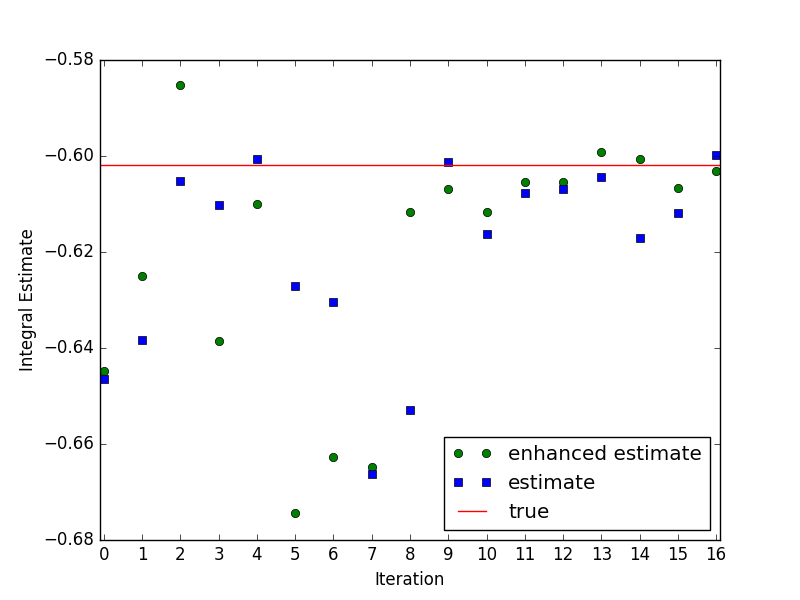}
\caption{Integral estimate at each iteration of the adaptive algorithm using the normal surrogate and the enhanced surrogate.}
\label{fig:ex1_int}
\end{figure}

\begin{figure}
\centering
\includegraphics[width=0.45\textwidth]{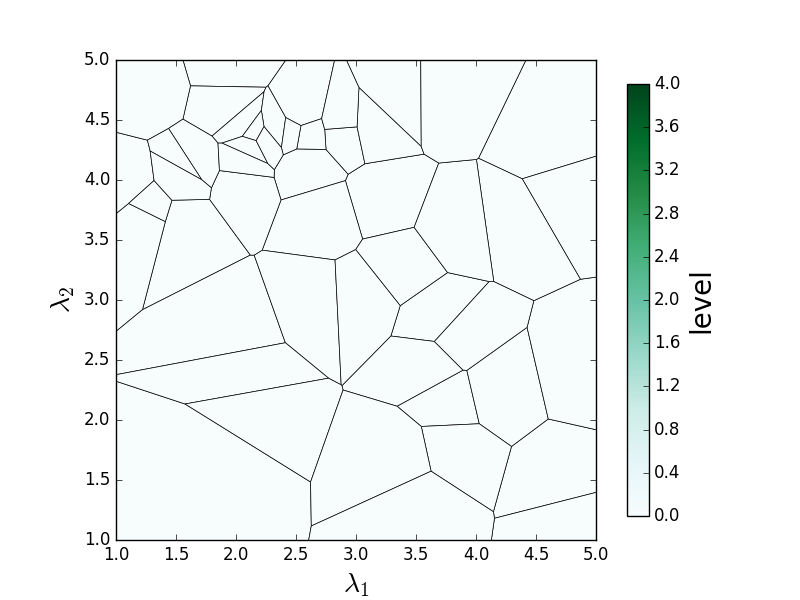}
\includegraphics[width=0.45\textwidth]{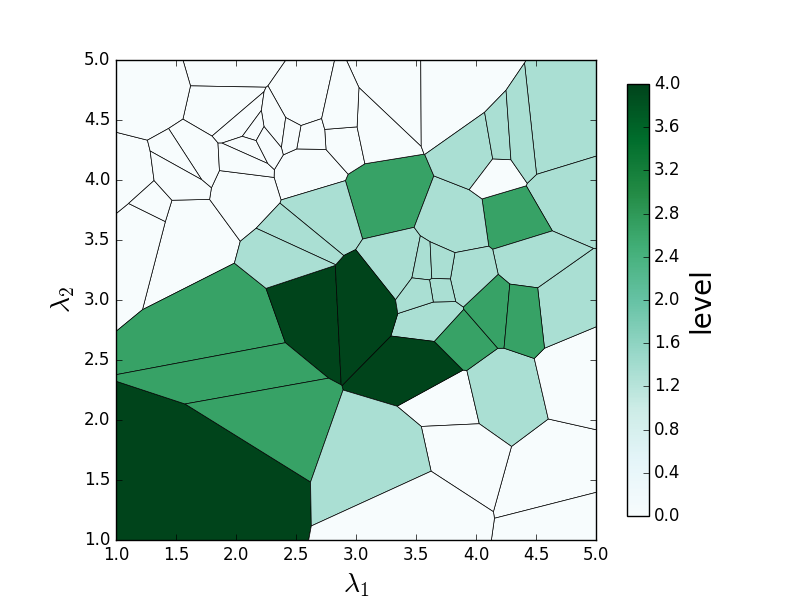}
\includegraphics[width=0.45\textwidth]{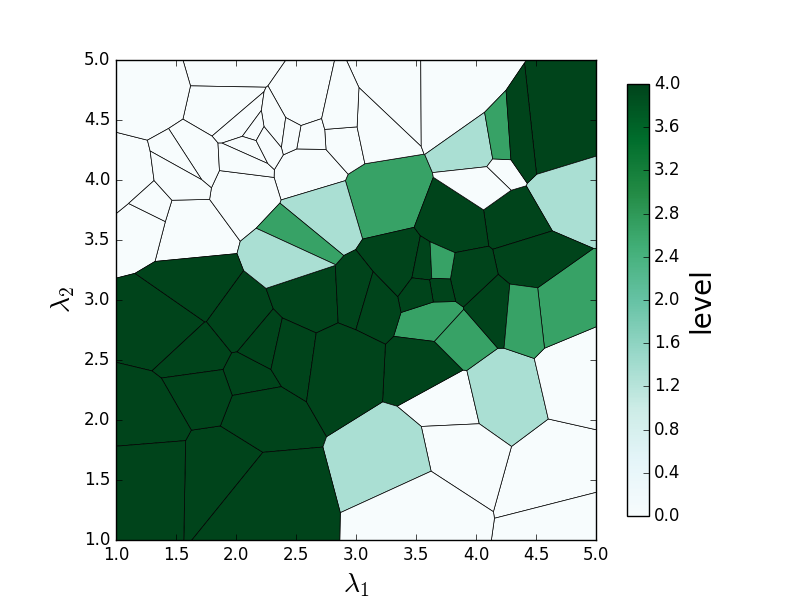}
\includegraphics[width=0.45\textwidth]{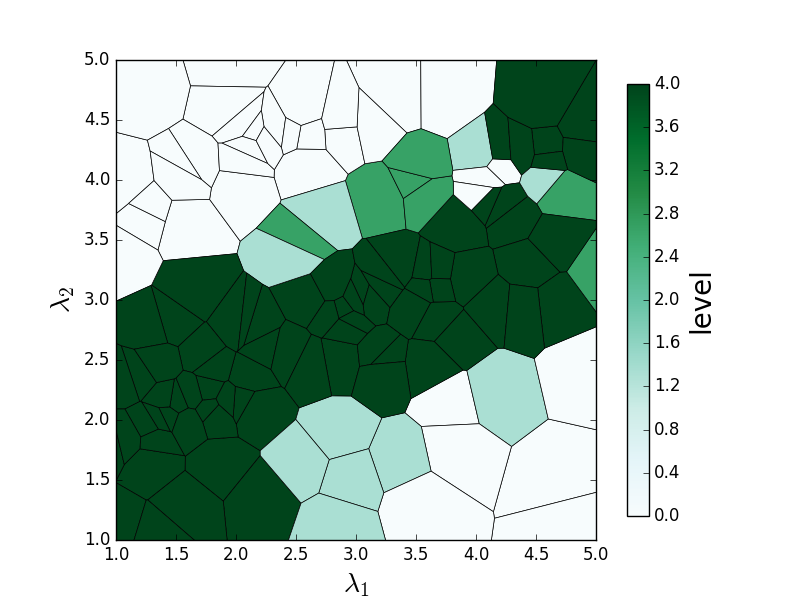}
\includegraphics[width=0.45\textwidth]{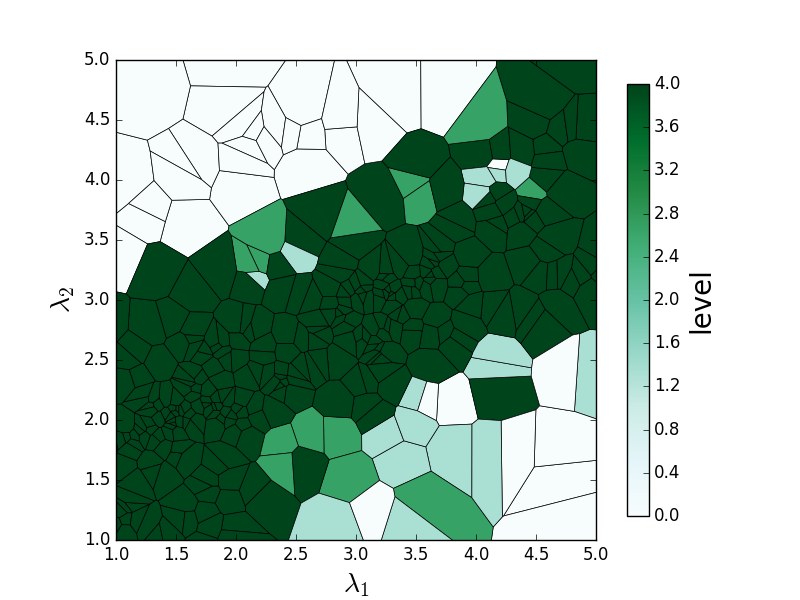}
\includegraphics[width=0.45\textwidth]{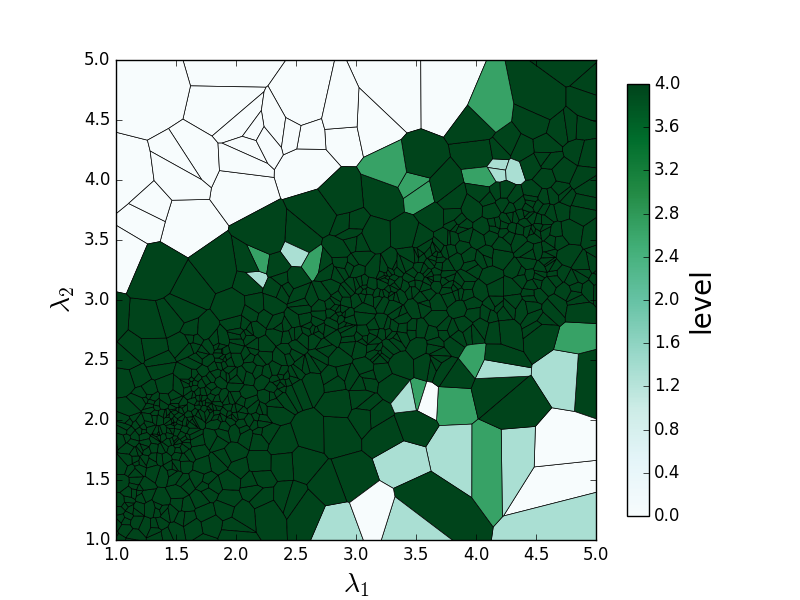}
\caption{Voronoi cells and model levels at iteration 0 (top left), 3 (top right), 6 (center left), 9 (center right), 12 (bottom left), and 15 (bottom right).}
\label{fig:ex1_levels}
\end{figure}

We initialize the surrogate with $N_0 =50$ uniform i.i.d. samples in $\pspace$.
The model is initially solved at level 1 with piecewise constants everywhere.
The adaptive algorithm (Alg. \ref{alg:ad_samp}) was used to adaptively update the surrogate and calculate the integral, with a relative tolerance of $\epsilon = 0.01$.
The algorithm was run until the convergence criterion was met.
The Monte Carlo method for computing error indicators for cheap numerical models was used.
Table \ref{tab:ex1_adapt_levels} shows the cumulative number of model evaluations at each level that were performed through each iteration for the first 16 iterations as well as the absolute error in the computation of the integral.
Figure \ref{fig:ex1_int} shows the calculated value of the integral at each iteration  calculated with the normal and enhanced surrogate.
Figure \ref{fig:ex1_levels} shows the Voronoi tessellation for the initial discretization and iterations 3, 6, 9, 12, and 15.

Notice that the first several iterations identify the cells of importance for the computation of the integral and mostly does level-refinement on these cells.
The calculated value of the integral changes a great deal at each iteration as the local deterministic error is removed from the surrogate.
Cancellations of error and errors not being removed at the same iteration at different locations of the domain cause the errors to oscillate for a few iterations.
By iteration 8, a small amount of h-refinement and much level-refinement has reduced the effect of much of the deterministic error by refining at the highest level in the regions that most influence the integral.
At this point, the integral is relatively well-approximated, and h-refinement in these important regions from iteration 8 onward causes the error to rapidly decrease.
The cells in the regions with little or no influence on the computation of the integral remain coarse and involve only solving the model at low levels, requiring little computational cost.
Notice that the error at iteration 12 (with 313 model solves at the highest level) gives a more accurate computation of the integral than with uniform refinement with 10000 model solves at the highest level.

This illustrates the performance of the adaptive surrogate construction.
Regions highly impacting the computation of the integral are slowly refined with a mixture of refinement methods.
Eventually, the effect of the deterministic error $\epsilon_h$ is mostly removed  and h-refinement takes over and the error rapidly decreases as the surrogate error $\epsilon_s$ has an increasingly smaller effect.
Our method avoids computationally expensive high-level model solves in areas with little influence on the integral and strategically performs these high-level model solves in areas with much influence on the computation of the integral.

\begin{table}
\centering
\begin{tabular}{|r||r|r|r|r|r|}
\hline
$N$ & Level 1 & Level 2 & Level 3 & Level 4 & Level 5 \\ \hline \hline
 10 & 1.15e-01 & 1.22e-01 & 1.21e-01 & 1.03e-01 & 9.90e-02\\ \hline 
 100 & 4.55e-02 & 1.80e-02 & 1.30e-02 & 1.20e-02 & 1.27e-02\\ \hline  
 1000 & 4.41e-02 & 9.59e-03 & 3.24e-03 & 2.32e-03 & 2.92e-03\\ \hline 
 10000 & 4.45e-02 & 9.37e-03 & 3.51e-03 & 2.62e-03 & 2.76e-03\\  
\hline
\end{tabular}
\caption{Average errors (over 20 runs) in the computed integrals using piecewise constant surrogates created with $N$ uniform generating points in $\pspace$ using the five levels of the mesh.}
\label{tab:ex1_const}
\end{table}

\begin{table}
\centering
\begin{tabular}{|r||r|r|r|r|r|}
\hline
$N$ & Level 1 & Level 2 & Level 3 & Level 4 & Level 5 \\ \hline \hline
10 &  1.56e-01 & 3.71e-01 & 2.27e-01 & 2.24e-01 & 9.26e-02 \\ \hline 
100 &  5.87e-02 & 2.44e-02 & 3.80e-02 & 2.07e-02 & 1.62e-02 \\ \hline 
1000 &  4.41e-02 & 8.48e-03 & 2.60e-03 & 2.50e-03 & 2.78e-03 \\ \hline 
 10000 & 4.38e-02 & 9.98e-03 & 2.96e-03 & 2.43e-03 & 2.27e-03 \\  
\hline
\end{tabular}
\caption{Average errors (over 20 runs) in the computed integrals using piecewise linear surrogates created with $N$ uniform generating points in $\pspace$ using the five levels of the mesh.}
\label{tab:ex1_lin}
\end{table}

\begin{table}
\centering
\begin{tabular}{|r||r|r|r|r|r|r|}
\hline
Iteration & Level 1 & Level 2 & Level 3 & Level 4 & Level 5 & Error \\
\hline \hline
0 & 50 &  0 &  0 &  0 &    0 & 4.30e-02 \\ \hline
1 & 54 & 11 &  0 &  0 &    0 &  2.32e-02\\ \hline
2 & 58 & 21 &  6 &  0 &    0 &  1.65e-02\\ \hline
3 & 58 & 32 & 12 &  4 &    0 &  3.68e-02 \\ \hline
4 & 58 & 32 & 23 &  9 &    2 & 8.22e-03 \\ \hline
5 & 59 & 32 & 25 & 20 &    6 & 7.25e-02\\ \hline
6 & 59 & 32 & 26 & 24 &   15 & 6.08e-02 \\ \hline
7 & 59 & 33 & 26 & 27 &   30 &  6.30e-02\\ \hline
8 & 60 & 35 & 31 & 34 &   55 & 9.88e-03\\ \hline
9 & 61 & 37 & 31 & 34 &   70 & 5.03e-03 \\ \hline
10 & 62 & 49 & 35 & 38 &  115 & 9.78e-03 \\ \hline
11 & 68 & 51 & 44 & 47 &  190 & 3.67e-03  \\ \hline
12 & 70 & 57 & 52 & 55 &  307 & 3.73e-03 \\ \hline
13 & 70 & 63 & 58 & 66 &  457 & 2.51e-03 \\ \hline
14 & 70 & 65 & 64 & 76 &  726 &  1.13e-03\\ \hline
15 & 70 & 65 & 67 & 79 &  854 & 4.90e-03 \\ \hline
16 & 70 & 73 & 80 & 95 & 1487 & 1.27e-03\\ \hline
%0 & 50 &  0 &  0 &  0 &   0 &  9.82e-02\\ \hline 
%1 & 53 & 16 &  0 &  0 &   0 &  6.54e-02\\ \hline 
%2 & 54 & 23 &  8 &  0 &   0 &  9.28e-02 \\ \hline 
%3 & 54 & 26 & 14 &  5 &   0 & 1.15e-01 \\ \hline 
%4 & 54 & 28 & 18 & 12 &   5 &  1.22e-01 \\ \hline 
%5 & 59 & 32 & 21 & 17 &   5 & 2.47e-02 \\ \hline 
%6 &  59 & 38 & 24 & 23 &  16 & 2.70e-02 \\ \hline 
%7 & 64 & 47 & 28 & 27 &  35 & 2.82e-02 \\ \hline 
%8 & 64 & 47 & 29 & 27 &  49 & 1.33e-02 \\\hline 
%9 & 64 & 48 & 32 & 28 &  75 &  9.87e-03\\ \hline 
%10 & 64 & 49 & 36 & 35 & 122 & 2.35e-03 \\ \hline 
%11 & 64 & 52 & 42 & 39 & 202 &  1.08e-02\\ \hline 
%12 & 64 & 53 & 46 & 45 & 313 &  1.81e-03  \\ 
\hline
\end{tabular}
\caption{Number of cumulative model evaluations at each level and absolute errors in the computation of the integral for 16 iterations of the adaptive algorithm for the 1D problem.}
\label{tab:ex1_adapt_levels}
\end{table}

%\begin{figure}
%\centering
%\includegraphics[width=0.8\textwidth]{Figures/hist_exact.png}
%\caption{Reference posterior distribution.}
%\label{fig:true_post}
%\end{figure}

\subsection{Elliptic PDE with Uncertain Conductivity}

\begin{table}
\centering
\begin{tabular}{|r||r|r|r|r|r|}
\hline
Iteration & Level 1 & Level 2 & Level 3 & Level 4 & Error \\ \hline \hline
   0 &  539 &   8 &   2 &   2 &  5.44e-02 \\ \hline
 100 &  696 &  33 &  10 &  12 & 2.40e-02\\   \hline
 150 &  764 &  49 &  16 &  22 & 2.29e-02 \\ \hline
 200 &  823 &  67 &  21 &  40 & 2.44e-02 \\ \hline
 250 &  884 &  85 &  28 &  54 & 2.48e-02 \\ \hline
 300 &  942 & 100 &  43 &  66 & 2.33e-02 \\ \hline
 350 & 1003 & 108 &  60 &  80 & 2.36e-02 \\ \hline
 400 & 1067 & 117 &  74 &  94 & 2.24e-02\\ \hline
 450 & 1126 & 126 &  93 & 106 &  2.18e-02\\ \hline
 500 & 1194 & 132 & 100 & 125 & 8.78e-03\\ \hline
 550 & 1272 & 148 & 104 & 131 & 2.88e-03  \\ \hline
 600 & 1345 & 169 & 111 & 140 &  2.96e-03 \\
%   0 &  500 &  0 &  0 &  0 &  1.13e-03 \\ \hline
%  40 &  572 &  5 &  2 &  2 & 6.97e-05\\ \hline
%  80 &  643 & 12 &  2 &  4 &  2.74e-04 \\ \hline
% 120 &  704 & 30 &  3 &  4 & 8.22e-05 \\ \hline
% 160 &  773 & 40 &  4 &  4 & 2.67e-04\\ \hline
% 200 &  848 & 44 &  5 &  4 & 6.16e-05\\ \hline
% 240 &  914 & 55 &  7 &  5 &  2.50e-04\\ \hline
% 280 &  965 & 60 &  9 & 27 & 2.62e-04\\ \hline
% 320 & 1020 & 77 & 11 & 33 & 3.25e-05 \\ \hline
\hline
\end{tabular}
\caption{Cumulative number of model evaluations at each level and absolute errors in the calculations of integrals for iterations of the adaptive scheme applied to the 2D PDE.}
\label{tab:ex2}
\end{table}

\begin{figure}
\centering
% Put back
\includegraphics[width=0.8\textwidth]{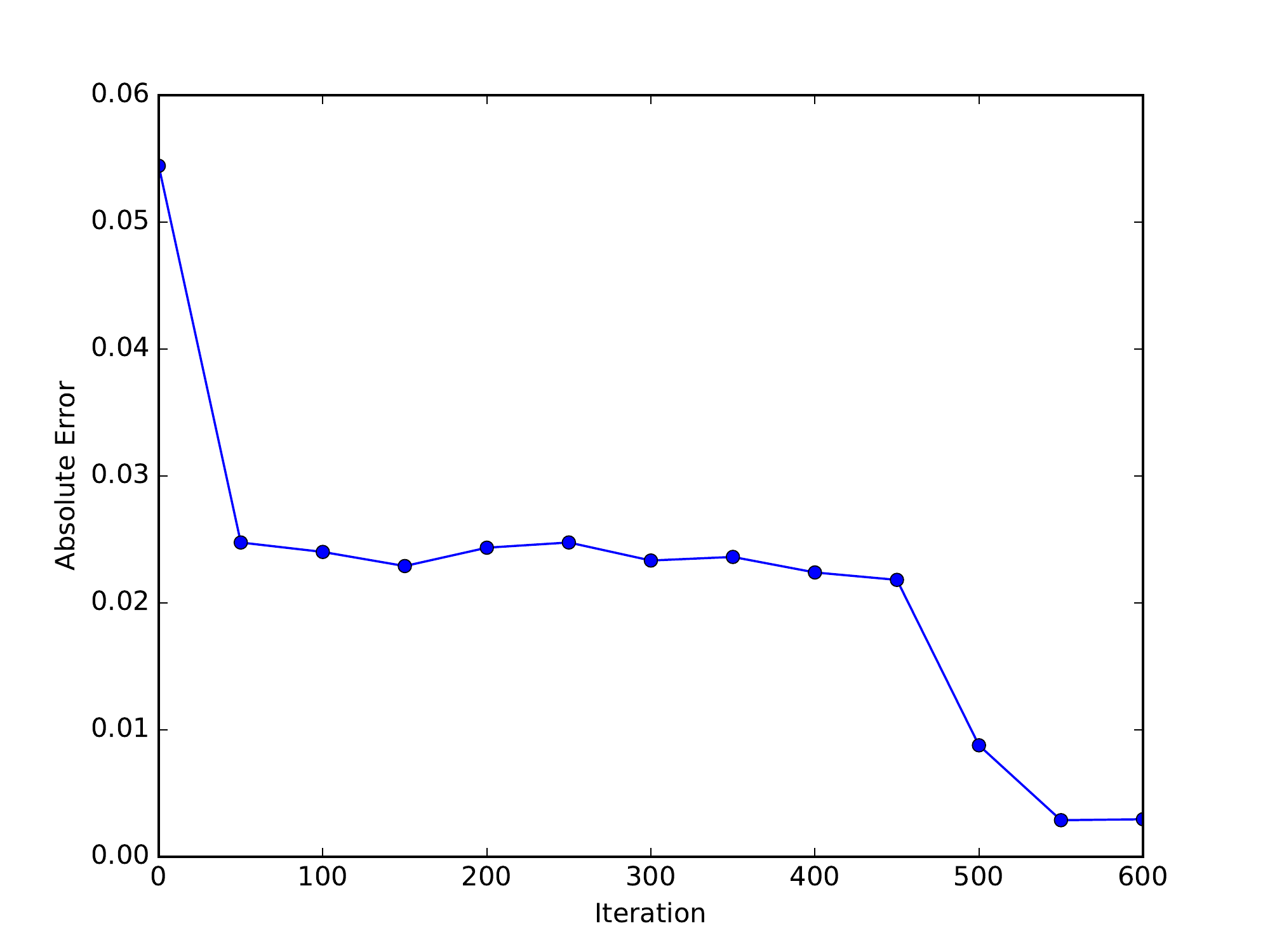}
\caption{Integral error in every 50th iteration for the 2D elliptic problem.}
\label{fig:ex2}
\end{figure}

We consider the elliptic boundary value problem on the unit square
\begin{eqnarray} \label{eq:ex2}
- \nabla \cdot \left( K(x,y) \nabla u(x,y) \right) &=& 0, \text{ for }  (x,y) \in (0,1)^2  \nonumber \\
u(0,y)  &=0&, \nonumber \\
 u(1,y) &=& 1, \nonumber \\
\nabla u(x,0) \cdot \mathbf{n} =\nabla u(x,1) \cdot \mathbf{n} &=& 0,
\end{eqnarray}
where $K(x,y)$ is a conductivity field
%Suppose $Y(x,y) = \log (K_s(x,y))$, is an instance of a random process, and let $Y'(x,y)$ be the mean-removed process. 
 which we treat as a random function.
$K$ belongs to an infinite-dimensional space, but truncating a Karhunen-Lo\`eve (K-L) Expansion is a classical option for deriving finite-dimensional parameterizations for $\log(K)$.
%Constructing the KLE first requires specification of a covariance function.
%This may be obtained, for instance, assuming a stationary random field and using a variogram on available data from a sparse set of boreholes.
We  construct the K-L Expansion of $Y(x,y)$ where $Y(x,y) := \log [K(x,y)]$.
Let $\bar{Y}(x,y)$ be the mean value of $Y(x)$, and that it has a exponential covariance $C$ with correlation lengths of 0.1 in both directions.
%Thus, it has the spectral decomposition
%\begin{equation}
%C(\boldsymbol{x_1},\boldsymbol{x_2}) = \sum_{n=0}^{\infty} \lambda_n f_n(\boldsymbol{x_1}) f_n(\boldsymbol{x_2})
%\end{equation}
%where $\lambda_n$ and $f_n(\boldsymbol{x})$ are the solutions to the homogeneous Fredholm integral equation of the second kind:
%\begin{equation}
%\int_{\boldsymbol{D}} C(\boldsymbol{x_1},\boldsymbol{x_2}) f_n(\boldsymbol{x_1}) d\boldsymbol{x_1} = \lambda_n f_n(\boldsymbol{x_2}).
%\end{equation}
%The eigenfunctions are orthogonal and form a complete set. 
Hence, $Y(x,y)$ can be written as
\begin{equation} \label{eq:KLE_infty}
Y(x,y) = \bar{Y}(x,y) + \sum_{n=0}^{\infty} \xi_n \sqrt{\lambda_n} f_n(x,y),
\end{equation}
where $\lambda_n$ and $f_n(x,y)$ are eigenpairs determined by $C$, and $\xi_n$ are standard normal random variables.
%The KLE of a Gaussian field has the further property that $\xi_n(\omega)$ are independent standard normal random variables \cite{le2010spectral}.
Truncating the series in Eq. \eqref{eq:KLE_infty} at the $N${th} term gives the finite-dimensional approximation
\begin{equation}\label{eq:KLE_truncate}
Y(x,y) \approx \bar{Y}(x,y) + \sum_{n=0}^{N} \xi_n \sqrt{\lambda_n} f_n(x,y).
\end{equation}
We use the first eight K-L terms (i.e. $N=8$) because the eigenvalues above this are observed to be negligible for this correlation length and take $\bar{Y}(x,y) = 0.05$.

Given K-L coefficients $\xi_i$, the system is discretized with the continuous Galerkin finite element method using linear elements on structured triangular grids using the open-source software FEniCS \cite{FEniCS_book, ans20553}. There are four levels of refinement of the mesh, which correspond to the model levels for the adaptive scheme, with 15x15, 21x21, 30x30, and 42x42 elements respectively.
The assembled linear systems are solved with a direct solver.
The uncertain parameters are the K-L coefficients $\lambda = [\xi_1, \xi_2,..., \xi_8]$, and thus the parameter 
space is $\pspace = \mathbb{R}^8.$
The stochastic inverse problem is defined with the forward map $Q: \pspace \rightarrow \mathbb{R}^4$, where  $Q(\lambda) = [Q_1(\lambda), Q_2(\lambda),Q_3(\lambda),Q_4(\lambda)]$, and $Q_1(\lambda) = u(0.25, 0.25)$, $Q_2(\lambda) = u(0.25 ,0.75)$, $Q_3(\lambda) = u(0.75, 0.25)$, $Q_4(\lambda) = u(0.75,0.75)$.
The corresponding adjoint problems are
\begin{eqnarray} \label{eq:ex2_adj}
- \nabla \cdot \left( K(x,y) \nabla \phi(x,y) \right) &=& \psi_i(x,y), \text{ for }  (x,y) \in [0,1]^2  \nonumber \\
\phi(0,y)  &=0&, \nonumber \\
 \phi(1,y) &=& 0, \nonumber \\
\nabla \phi(x,0) \cdot \mathbf{n} =\nabla \phi(x,1) \cdot \mathbf{n} &=& 0,
\end{eqnarray}
where $\psi_i(x,y)$ are steep Gaussians approximating a Dirac delta at the evaluation points of $Q_i$.
The adjoint problems are solved on the same finite element meshes, but with quadratic elements.
Using an enriched space for the adjoint problem is required with finite elements for performing error estimates.
Error estimates and gradients can be calculated using Equations \ref{eq:pls_QoI_error} and \ref{eq:pls_QoI_deriv}.

We pose the stochastic inverse problem as a Bayesian inverse problem.
The data is $y_{data} = [0.2803, 0.2693, 0.8114, 0.6506]$.
We assume a standard normal prior on $\pspace$ and mean-zero Gaussian noise $\eta \backsim \mathcal{N}(0,0.0025)$ in each component.
The function $f$ that we want to integrate is 
$$f(\lambda) = \int_{0.4}^{0.6} \int_{0.4}^{0.6} K(x,y, \lambda) dx dy$$
 the upscaled (through volume averaging) conductivity over the block $[0.4, 0.6]^2$.
The method calculating error indicators that is used is the method for ``expensive" models.
The problem is solved with the Metropolis-Hastings algorithm.
A reference solution using $10^5$ model evaluations at the highest level is used for comparison with results from the adaptive scheme.

The adaptive scheme was initialized with 500 i.i.d. samples in $\pspace$ with respect to the standard normal distribution.
The model was solved at the lowest level for the initial discretization.
The adaptive scheme was run until it met the stopping criterion with $\epsilon = 10^{-4}$.
The adaptive algorithm terminated after 600 iterations with an absolute error of $2.96 \times 10^{-3}$.
Table \ref{tab:ex2} shows the number of average number of cumulative model evaluations at each level for every 50 iterations as well as the average absolute error in the calculation of the integral.
We see that in the early iterations, much h-refinement and a small amount of level-refinement occurs.
In the moderately high dimensional space $\pspace$, this is adding new generating points in the areas where the posterior has larger values.
The error quickly decreases as samples are added in these regions.
In the middle steps of the process, the error fluctuates as the surrogate is refined, sometimes with cancellations of error occurring.
After around 450 iterations, h-refinement at the lowest levels continue and level-refinement increases in the areas where deterministic model error has a larger effect and a great reduction in error occurs.
The process terminates with 1765  model evaluations, with the grand majority at the lowest level.
This indicates that in this problem, error due to the surrogate was generally more important than the deterministic error in the calculation of the QoI.
This result is likely different from the previous example because of the dimension of $\pspace$.
Using a surrogate in higher dimensions can introduce more error than in lower ones.
h-refinement at the lowest level was mostly able resolve the surrogate enough for an accurate solution.
Expensive level 4 evaluations only had to be done 140 times, in the areas where the deterministic model error greatly affects the computation of the integral.
The adaptive scheme was able to provide a highly accurate computation of the integral with a relatively small number of model evaluations.
Most of the model evaluations could be done with the cheap low-level model.
This illustrates the methods ability to make predictions under uncertainty accurately with much lower computational costs than a non-adaptive method.

\subsection{Predator-Prey Model}

\begin{figure}
\centering
% Put back
\includegraphics[width=0.8\textwidth]{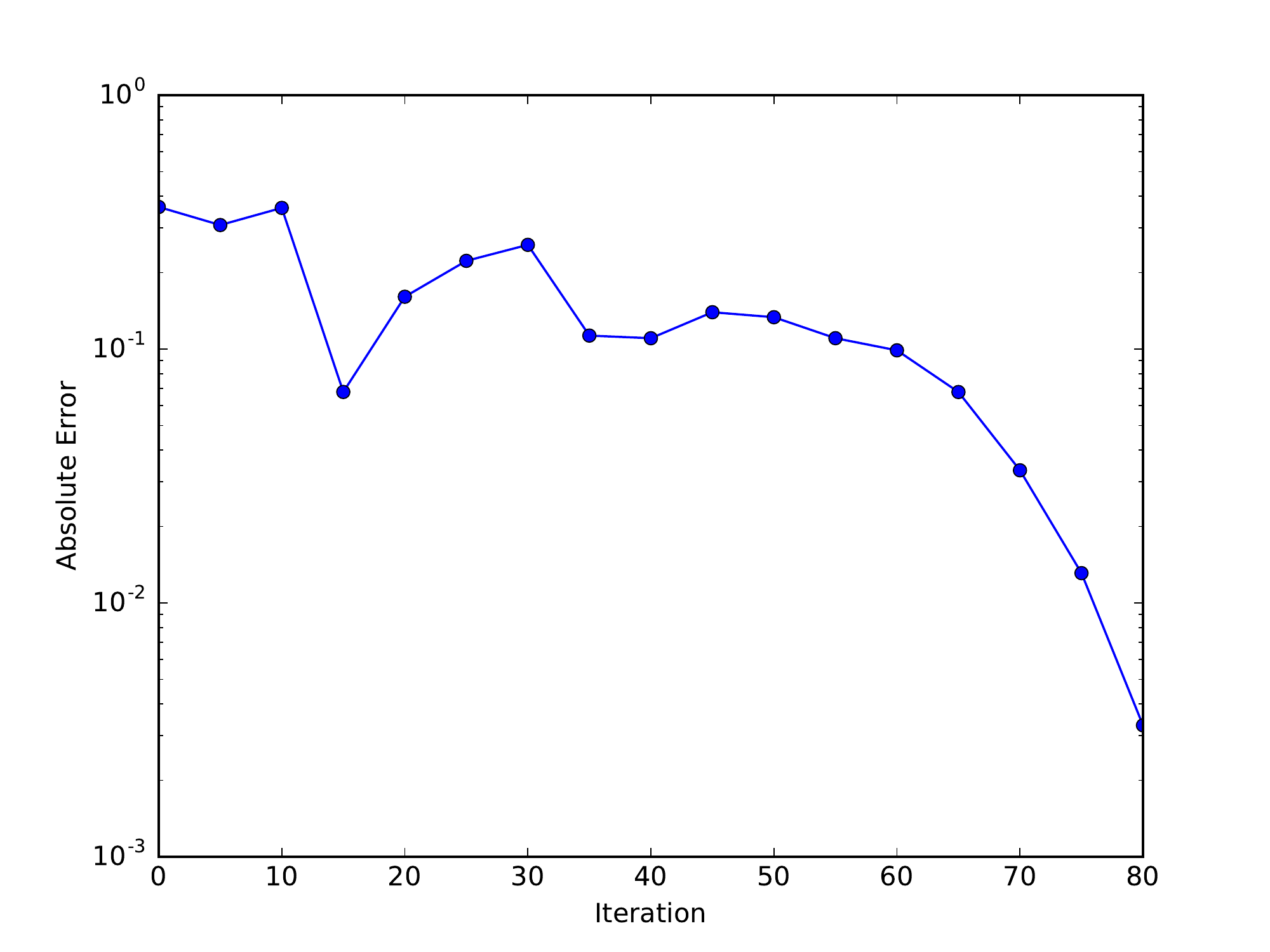}
\caption{Integral error in every 5th iteration for the predator-prey model.}
\label{fig:ex3}

\end{figure}

\begin{table}
\centering
\begin{tabular}{|r||r|r|r|r|r|r|}
\hline
Iteration & Level 1 & Level 2 & Level 3 & Level 4 & Error & Run Time (s) \\ \hline \hline
  0 & 100 &  0 &  0 &    0 &  3.63e-01 & 1.16e+00\\ \hline
  5 & 103 & 15 & 12 &    7 & 3.08e-01  & 1.82e+01 \\ \hline
 10 & 103 & 20 & 16 &   23 & 3.60e-01 &  5.13e+01 \\ \hline
 15 & 104 & 23 & 18 &   45 &  6.78e-02 & 9.60e+01 \\ \hline
 20 & 104 & 26 & 19 &   79 & 1.61e-01 & 1.65e+02 \\ \hline
 25 & 105 & 26 & 20 &  114 &  2.23e-01 &  2.35e+02\\ \hline
 30 & 105 & 27 & 20 &  159 & 2.57e-01  &  3.26e+02 \\ \hline
 35 & 105 & 27 & 20 &  245 & 1.13e-01 & 4.99e+02\\ \hline
 40 & 105 & 27 & 20 &  385 & 1.10e-01 &  7.80e+02 \\ \hline
 45 & 105 & 27 & 20 &  599 &  1.40e-01 & 1.21e+03  \\ \hline
 50 & 105 & 27 & 20 &  833 & 1.33e-01 & 1.68e+03 \\ \hline
 55 & 106 & 27 & 20 & 1098 & 1.10e-01  & 2.21e+03 \\ \hline
 60 & 106 & 27 & 20 & 1436 & 9.89e-02  & 2.89e+03  \\ \hline
 65 & 106 & 27 & 20 & 2128 & 6.77e-02 & 4.28e+03 \\ \hline
 70 & 106 & 27 & 20 & 3322 & 3.33e-02  & 6.68e+03  \\ \hline
 75 & 106 & 27 & 20 & 5472 & 1.31e-02  & 1.10e+04 \\ \hline
 80 & 106 & 27 & 20 & 6875 & 3.29e-03 & 1.38e+04  \\ \hline
\end{tabular}
\caption{Cumulative number of model evaluations at each level and absolute errors in the calculations of integrals for iterations of the adaptive scheme applied to the predatory-prey model.}
\label{tab:ex3}
\end{table}

\begin{table}
\centering
\begin{tabular}{|r|r|r|}
\hline
Num. of Evals. & Error & Run Time (s) \\ \hline
100 &  2.56e-01 & 2.01e+02 \\ \hline
1000 &  2.53e-01 & 2.01e+03 \\ \hline
10000 &  2.32e-01 & 2.01e+04 \\ \hline
\end{tabular}
\caption{Absolute errors and run times for the predator-prey model with uniform refinement at the finest level.}
\label{tab:ex3:2}
\end{table}

We consider the classical predator-prey model, a nonlinear system of ordinary differential equations defined by the Lotka-Volterra Equations
\begin{equation} \label{eq:predprey}
{\begin{aligned}{\frac {dx}{dt}}&=\alpha x-\beta xy\\[6pt]{\frac {dy}{dt}}&=\delta xy-\gamma y, \end{aligned}}
\end{equation}
for $t \in [0,T]$ and with initial conditions $x(0)=x_0$ and $y(0)=y_0$.
$x$ and $y$ represent the population of species of prey and predators respectively at time $t$, and the other parameters describe the population dynamics.
The six parameters $\alpha$, $\beta$, $\delta$, $\gamma$, $x_0$, and $y_0$ are unknown and uncertain.
For solving the model, the backward Euler method is used with Newton's method being used to solve the nonlinear system at each step, ensuring stability even with large time steps.
The calculated solutions are $x_h$ and $y_h$.
The QoIs are the populations of both species at times $t=5$ and $t=10$.
The adjoint problems associated with these QoI involve solving the linearized system
\begin{equation}
\frac{d\boldsymbol{\phi}}{dt} = J(x_h(t),y_h(t))\boldsymbol{\phi}
\end{equation}
backward in time for $t \in [T,0]$ where $J$ is the Jacobian of the RHS of Equation \ref{eq:predprey}, with $T$ and initial conditions corresponding to the respective QoI.
We solve the adjoint system using the Crank-Nicolson method for time integration.
Note that each time step of the adjoint problem only involves solving a 2x2 linear system.
Error estimates and gradients can be calculated using Equations \ref{eq:pls_QoI_error} and \ref{eq:pls_QoI_deriv}, where the inner product $\langle \cdot, \cdot \rangle$, is the space-time inner product $\langle \mathbf{z}_1, \mathbf{z}_2 \rangle = \int_0^T \mathbf{z}_2^T \mathbf{z}_1 dt$.
The time integrals are approximated with the midpoint rule over each time step.
We consider four levels of the model with different time steps $\Delta t = 0.25, 0.1, 0.01,$ and  $0.001$ respectively.

The parameter domain is $\pspace = [1,2]^6$ with $\lambda \in \pspace$ defined by $\lambda = [\alpha, \beta, \delta, \gamma, x_0, y_0]$.
The stochastic inverse problem is a Bayesian inverse problem with the map $Q: \pspace \rightarrow \mathbb{R}^4$,
 where $Q(\lambda) = [Q_1(\lambda), Q_2(\lambda),Q_3(\lambda),Q_4(\lambda)]$, and $Q_1(\lambda) = x(5)$, $Q_2(\lambda) = y(5)$, $Q_3(\lambda) = x(10)$, $Q_4(\lambda) = y(10)$.
The data is $y_{data} = [1, 1.8, 0.5, 1.4]$.
We assume a uniform prior on $\pspace$ and mean-zero Gaussian noise $\eta \backsim \mathcal{N}(0,0.065)$ in each direction.
The function that we want to integrate with respect to the posterior is $f(\lambda) = x_0/y_0$, the ratio of the initial populations.

The adaptive scheme was initialized with 100 samples chosen uniformly in $\pspace$.
The model was solved at level 1 at each sample, and the iterative process was run until it converged with $\epsilon = 10^{-3}$, which took 80 iterations.
The errors are calculated with a reference value coming from the solution to the Bayesian inverse problem with $10^5$ samples at the highest model level (with a corresponding run time of $2 \times 10^6$ s).
Table \ref{tab:ex3} shows the number of cumulative model evaluations at each level for every 5th iteration as well as the absolute error in the calculation of the integral.
Figure \ref{fig:ex3} shows the error at each of these iterations.
As with the first example, notice that for the first 50 iterations the error oscillates without reducing substantially.
During these iterations, a small amount of h-refinement and a large amount of level-refinement is being performed, slowly identifying the regions in $\pspace$ with significant impact on the computation of the integral.
These refinements are done with relatively few model evaluations at any level.
After approximately 50 iterations, the effect of the model discretization error $\epsilon_h$ on the computation of the integral has mostly been removed.
The dominating error contribution is from the local inaccuracy of the surrogate model, so the adaptive scheme continues with h-refinement at the highest model level in the important regions.
The error in the integral approximations steadily decreases as the surrogate error is reduced with the h-refinement.
The method terminates at the 80th iteration with approximately 7000 model evaluations.
Table \ref{tab:ex3:2} shows errors and run times for uniform refinement using the highest level model.
Using 10,000 fine model evaluations at a much higher computational cost than the adaptive method results in a much greater error (almost two orders of magnitude).
This illustrates the computational savings of the adaptive method.

\subsection{A Vibroacoustics Application}
\begin{figure}
\centering
% Put back
\includegraphics[width=0.5\textwidth]{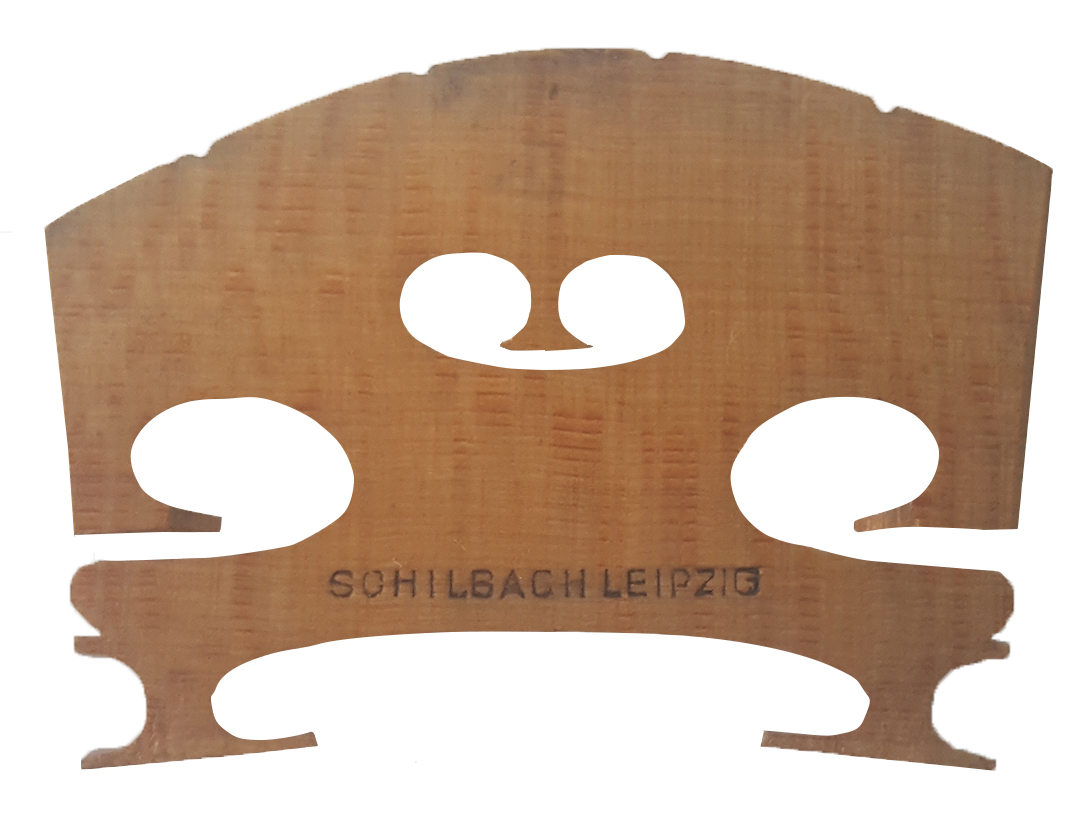}
\caption{A standard violin bridge.}
\label{fig:e4}
\end{figure}

Vibroacoustic applications often involve complex domains such as bridges, mechanical devices, stereo equipment, and musical instruments.
We consider the bridge of a violin (shown in Figure \ref{fig:e4}), which plays a critical role in transmitting the vibration of strings to the body of the instrument.
Violin bridges have complex curved geometries that are difficult to mesh with classical methods.
Because of this, we utilize an isogemetric mortar discretization.
Also, for computational efficiency we use a reduced basis method with saddle point approach which allows efficient construction of the linear systems given material parameters.
A full description and analysis of the reduced basis isogeometric mortar discretization is done by Horger et al. \cite{horger2017}.

On the three-dimensional violin bridge we consider a problem of linear elasticity
\begin{equation}
-div (\sigma(u))= f,
\end{equation}
where the strain $\sigma(u)$ depends on the material laws of the bridge and f are body forces.
The stress-strain relationship by Hooke's law $\sigma(u) = \mathbb{C} \epsilon(u)$, where $\epsilon(u) = \left( \nabla u + \nabla u^T \right)/2$.
The material is orthotropic, and the stiffness tensor is given by
$$\mathbb{C} = \left[ \begin {array}{cccccc} A_{11}&A_{12}&A_{13}&0&0&0\\ \noalign{\smallskip}A_{21}&A_{22}&A_{23}&0&0&0\\ \noalign{\smallskip}
A_{31}&A_{32}&A_{33}&0&0&0\\ \noalign{\smallskip}
0&0&0&G_{yz}&0&0\\ \noalign{\smallskip}
0&0&0&0&G_{zx}&0\\ \noalign{\smallskip}
0&0&0&0&0&G_{xy}
\end {array}
 \right],$$
where $G_{yz}$, $G_{zx}$, and $G_{xy}$ are the shear moduli, and the entries $A_{ij}$ depend on the elastic moduli $E_x$, $E_y$, $E_z$ and the Poisson's ratios $\nu_{xy}$, $\nu_{yz}$, $\nu_{zx}$.
The exact formulation of $\mathbb{C}$ is shown in \cite{rand2007}.
Assuming a known thickness and mass density of the bridge, there are thus nine material parameters for the structure:
$\lambda = [E_x, E_y, E_z, \nu_{xy}, \nu_{yz}, \nu_{zx}, G_{yz}, G_{zx}, G_{xy}]$.
We assume that the bridge is made out of \textit{Fagus sylvatica}, the common beech.
There are known ranges for the material parameters which are shown in Table \ref{tab:violin}.
Let the domain of possible parameters $\pspace$ be the tensor product of these intervals.

\begin{table}
%\centering
\small
\begin{tabular}{|r|r|r|r|r|r|r|r|r|r|} 
\hline
 & $E_x$ [MPa] & $E_y$ [MPa] & $E_z$ [MPa] & $\nu_{xy}$& $\nu_{yz}$& $\nu_{zx}$& $G_{yz}$ [MPa] & $G_{zx}$ [MPa]  & $G_{xy}$ [MPa]  \\ \hline \hline
min & 13,000 & 1,500 & 750 & 0.3 & 0.03 & 0.4 & 100 & 500 & 1000 \\ \hline 
max & 15,000 & 3,000 & 1,500 & 0.4 & 0.06 & 0.5 & 1,000 & 1,500& 2,000 \\ \hline
\end{tabular}
\caption{Ranges of material parameters for the violin bridge.}
\label{tab:violin}
\end{table}

\begin{table}
\centering
\begin{tabular}{|r|r|r|r|r|r|} \hline
Iter. & Level 1 & Level 2 & Level 3 & Rel. Error & Run Time (s)\\ \hline \hline
%0 & 100 & 0 & 0 & 0 & 5.48e-03 & 6,010 \\ \hline
0 & 100 & 0 & 0 & 1.04e-02 & 6,544 \\ \hline
1 & 102 & 70 & 0 & 5.69e-03 & 20,572\\ \hline
2 & 103 & 73 & 49 & 4.70e-04 & 35,708 \\ \hline
\end{tabular}
\caption{Results from two iterations of the adaptive scheme for the violin bridge problem.} \label{tab:violin2}
\end{table}

\begin{table}
\centering
\begin{tabular}{|r|r|r|} \hline
Num. of Evals. & Rel. Error & Run Time (s)  \\ \hline \hline
10 &  2.93e-02 &  601 \\ \hline
100 & 3.45e-02 & 6,010 \\ \hline
1000 &1.37e-02  &  60,100 \\ \hline
\end{tabular}
\caption{Results with uniform refinement at level 1 for the violin bridge problem.} \label{tab:violin3} 
\end{table}

We assume a known homogeneous force on the top of the bridge from the strings, homogeneous Dirichlet conditions where the bridge attaches to the body, and Neumann conditions on the remaining boundaries.
The reduced basis isogeometric mortar discretization results in a large sparse linear saddle point system.
We consider three mesh levels and a maximally smooth $p=3$ discretization and a $p=4$ discretization containing the $p=3$ space on each of these meshes.
The adjoint problems, error representations, and derivative calculations are constructed by the methods described in Section \ref{sec:adj}.
The forward problems are solved with the $p=3$ discretizations and the adjoint problems are usually solved with the $p=4$ discretizations.
For the coarsest mesh, the $p=3$ system has 9,132 degrees of freedom (DOF) and the $p=4$ system has 22,635 DOF.
For the middle mesh, the $p=3$ system has 15,468 DOF and the $p=4$ system has 47,985 DOF.
For the finest mesh, we only consider the $p=3$ system which has 47,985 DOFs.
Because the solution is already so accurate, we do not calculate error estimates and solve the adjoint problems with the same matrix to calculate derivatives.

The QoI map $Q$ has five components.
$Q_1$ and $Q_2$ are the average x-displacements on the left and right boundaries of the bridge respectively.
$Q_3$, $Q_4$ and $Q_5$ are the average displacements on the front face of the bridge in the x, y, and z directions respectively.
This QoI map is used to define a Bayesian inverse problem with artificial data generated by solving the fine model with reference parameters and adding noise, $y = [1.79 \times 10^{-4}, -6.57\times 10^{-4}, -4.17\times 10^{-3}, -1.31\times 10^{-3}, -6.61\times 10^{-2}]$, where the displacement data is measured in cm.
The solution is a posterior measure $P_{\pspace}$ on $\pspace$.
We assume a uniform prior on $\pspace$ and mean-zero Gaussian noise at a level of 20\%.
The goal is to predict the expected average vertical displacement on the top of the bridge, i.e. $f(x,y,z, \lambda) = \int_{\Omega_{top}} u_z(x,y,z, \lambda) ds$.

Our adaptive algorithm was initialized with 100 uniform i.i.d. samples in $\pspace$ with a piecewise linear surrogate ($p_i=1$).
%Because of the difficulty and computational expense of evaluating $f$, the error indicators $\widehat{E}_{int,i}$ are calculated using the method for expensive models.
The model is computationally expensive, with one full evaluation including adjoint solves taking 67.7 s, for level 1, 198.5 s for level 2, and 295.4 s for level 3 approximately, so the error indicators $\widehat{E}_{int,i}$ are calculated using the method for expensive models.
A reference computation of the integral of -0.01412 cm was calculated using an error-corrected piecewise linear surrogate generated from 1000 uniform samples solved at level 2, which relates to a run time of 55.1 hours.
The adaptive scheme was run for two iterations before converging.
Table \ref{tab:violin2} shows results from each iteration including the cumulative number of model evaluations, the relative error compared to the reference, and the cumulative run time.
With the initial surrogate, the relative error in the calculation of the integral is already quite small.
This small error is due to the fact that for this problem, the QoI response with respect the the parameters are close to linear locally, which is not evident a priori.
Hence, the piecewise linear surrogate does not introduce much error into the prediction directly.
Most of the error is because the error in the numerical computation of the QoI with the coarse discretization.
Between iterations 0 and 1, level-refinement is performed on 70 cells to reduce the effect of the discretization error, and a small amount of h-refinement refines the Voronoi tessellation.
These refinements cause a large reduction in the error in the predicted value.
Between iterations 1 and 2, more level-refinement is performed on 49 cells requiring the fine model to be solved, and a small amount of h-refinement occurs at levels 1 and 2.
These refinements decrease the effect of deterministic model errors even more, and the relative error decreases by an order of magnitude, and the algorithm terminates.

The piecewise polynomial surrogate, combined with our adaptive algorithm have proven to be highly successful for this problem.
The piecewise linear surrogate does a good job of approximating the QoI response surface with is locally close to linear.
The adaptive scheme identifies the areas where deterministic model error is polluting the prediction and refines accordingly.
In a relatively small computation time (approximately 10 hours of CPU time), a highly accurate calculation of the predicted value.
Table \ref{tab:violin3} shows corresponding errors and run times for uniform h-refinement of the surrogate using model level 1.
The error does slightly decrease as the the number of model evaluations increases, but at a slow rate.
The effect of deterministic model error is not being reduced as it is with the adaptive method.
Even with 1000 model evaluations and almost double the CPU time as with the adaptive method, the error in the prediction is still two orders of magnitude higher than with the adaptive scheme.

%%%%%%%%%%%%%%%%%%%%%%%%%%%%%%%%%%%%%%%%%%%%%%%%%%%%%%%%%%%%%%%%%%%%
%%%%%%%%%%%%%%%%%%%%%%%%%%%%%%%%%%%%%%%%%%%%%%%%%%%%%%%%%%%%%%%%%%%%

\section{Conclusions}\label{sec:Conclusions}
We have presented a method for goal-oriented adaptive surrogate construction for prediction under uncertainty.
A general class of surrogate models for response surfaces based on piecewise polynomial approximations on Voronoi tessellations forms the basis of this adaptive strategy.
The solution of adjoint problems is used to enhance these surrogates via derivative information which is used to increase the local polynomial order of the approximation and via a posteriori error estimates for QoIs.
These enhancements are used to create two levels of surrogates from which local error indicators are derived.
Computational algorithms for estimating these error indicators are also presented.
The error indicators are used to guide p-refinement, level-refinement, and h-refinement of the surrogates.
Such refinements improve both the regular and enhanced surrogates.
The surrogates and refinement strategies are combined in an iterative method for surrogate construction which reduces the effect of various types of discretization and surrogate errors on the computations of integrals corresponding to predictions under uncertainty.

The presented method is applied to four example problems of varying complexity; however, even for relatively simple forward models, the map between parameters and QoIs is often highly nonlinear and quite complex.
The results show that the method is successful in accurate computations of the integrals of interest with a relatively cheap computational cost.
It is important to note that the method was tested on problems of moderate parameter dimension.
Certain attributes of the algorithm are not tenable for very high-dimensional problems (the ``curse of dimensionality"), which is certainly a drawback.
However, there has been much work recently in techniques for effectively reducing the dimension of stochastic inverse problems (e.g. active subspaces \cite{CDW_ActiveSubspace_2014} and reduced basis methods \cite{le2010}) which can potentially be applied to higher dimensional problems to reduce the effective dimension to a moderate values for which the presented method is feasible.
The combination of the method presented here with such dimension reduction techniques is left to future work.

\section*{Acknowledgments}
Financial support was provided by the the German Research Foundation (DFG, Project WO 671/11-1).
The authors acknowledge Linus Wunderlich for his help in setting up the violin bridge problem and for the violin bridge photograph.

\bibliographystyle{siam}
\bibliography{references}

\end{document}